\documentclass[11pt]{article}
\usepackage{latexsym, epsfig, amssymb, amsmath, amsthm, graphicx, mathrsfs}
\usepackage[dvips]{color}
\usepackage{times}
\usepackage{anysize}

\marginsize{1 in}{1 in}{0.35 in}{1.2 in} %
\makeatletter

\newtheorem{lemma}{Lemma}
\newtheorem{theorem}{Theorem}

\makeatother

\newcommand{\be}{\mbox{\bf e}}

\newcommand{\br}{\mbox{\bf r}}
\newcommand{\bu}{\mbox{\bf u}}
\newcommand{\bv}{\mbox{\bf v}}

\newcommand{\bx}{\mbox{\bf x}}
\newcommand{\by}{\mbox{\bf y}}
\newcommand{\bz}{\mbox{\bf z}}
\newcommand{\bA}{\mbox{\bf A}}

\newcommand{\bD}{\mbox{\bf D}}

\newcommand{\bR}{\mbox{\bf R}}
\newcommand{\bS}{\mbox{\bf S}}
\newcommand{\bU}{\mbox{\bf U}}
\newcommand{\bV}{\mbox{\bf V}}
\newcommand{\bW}{\mbox{\bf W}}
\newcommand{\bX}{\mbox{\bf X}}

\newcommand{\bZ}{\mbox{\bf Z}}
\newcommand{\bone}{\mbox{\bf 1}}
\newcommand{\bzero}{\mbox{\bf 0}}
\newcommand{\bveps}{\mbox{\boldmath $\varepsilon$}}
\newcommand{\bbeta}{\mbox{\boldmath $\beta$}}

\newcommand{\bxi}{\mbox{\boldmath $\xi$}}
\newcommand{\bleta}{\mbox{\boldmath $\eta$}}

\newcommand{\bzeta}{\mbox{\boldmath $\zeta$}}
\newcommand{\bomega}{\mbox{\boldmath $\omega$}}

\newcommand{\hbbeta}{\widehat\bbeta}
\newcommand{\hbeta}{\widehat\beta}

\newcommand{\heq}{\ \widehat=\ }

\newcommand{\var}{\mathrm{var}}
\newcommand{\cov}{\mathrm{cov}}

\newcommand{\Sig}{\mathbf{\Sigma}}
\newcommand{\veps}{\varepsilon}

\newcommand{\diag}{\mathrm{diag}}

\newcommand{\SIS}{{ SIS}}
\newcommand{\ISIS}{{ ISIS}}
\newcommand{\Lasso}{{Lasso}}

\def\t{^T}

\def\toas{\overset{\mathrm{a.s.}}{\longrightarrow}}
\def\deq{\overset{\mathrm{(d)}}{=\hspace{-0.02 in}=}}

\begin{document}

\title{\textbf{Sure Independence Screening for Ultra-High Dimensional Feature Space} %
\thanks{Financial support from the NSF grants DMS-0354223, DMS-0704337 and DMS-0714554,
and the NIH grant R01-GM072611 is gratefully acknowledged. We are
grateful to the anonymous referees for their constructive and
helpful comments. Address for correspondence: Jianqing Fan,
Department of Operations Research and Financial Engineering,
Princeton University, Princeton, NJ 08544. Phone: (609) 258-7924.
E-mail: jqfan@princeton.edu.}
\date{\today}
\author{Jianqing Fan\\
\small Department of Operations Research and Financial Engineering\\
\small Princeton University\\
\and Jinchi Lv\\
\small Information and Operations Management Department\\
\small Marshall School of Business\\
\small University of Southern California\\
} %
}

\maketitle

\vspace*{-0.5 in}

\begin{abstract}

Variable selection plays an important role in high dimensional
statistical modeling which nowadays appears in many areas and is key
to various scientific discoveries. For problems of large scale or
dimensionality $p$, estimation accuracy and computational cost are
two top concerns. In a recent paper, Candes and Tao (2007) propose
the Dantzig selector using $L_1$ regularization and show that it
achieves the ideal risk up to a logarithmic factor $\log p$. Their
innovative procedure and remarkable result are challenged when the
dimensionality is ultra high as the factor $\log p$ can be large and
their uniform uncertainty principle can fail.

Motivated by these concerns, we introduce the concept of sure
screening and propose a sure screening method based on a correlation
learning, called the Sure Independence Screening (SIS), to reduce
dimensionality from high to a moderate scale that is below sample
size. In a fairly general asymptotic framework, the correlation
learning is shown to have the sure screening property for even
exponentially growing dimensionality. As a methodological extension,
an iterative SIS (ISIS) is also proposed to enhance its finite
sample performance. With dimension reduced accurately from high to
below sample size, variable selection can be improved on both speed
and accuracy, and can then be accomplished by a well-developed
method such as the SCAD, Dantzig selector, Lasso, or adaptive Lasso.
The connections of these penalized least-squares methods are also
elucidated.
\end{abstract}

\textit{Short title}: Sure Independence Screening

\newpage
\section{Introduction}

\subsection{Background}
Consider the problem of estimating a $p$-vector of parameters
$\bbeta$ from the linear model
\begin{equation} \label{011}
\by=\bX\bbeta+\bveps,
\end{equation}
where $\by=\left(Y_1,\cdots,Y_n\right)\t$ is an $n$-vector of
responses, $\bX=\left(\bx_1,\cdots,\bx_n\right)\t$ is an $n\times p$
random design matrix with i.i.d. $\bx_1,\cdots,\bx_n$,
$\bbeta=\left(\beta_1,\cdots,\beta_p\right)\t$ is a $p$-vector of
parameters, and $\bveps=\left(\veps_1,\cdots,\veps_n\right)\t$ is an
$n$-vector of i.i.d. random errors. When dimension $p$ is high, it
is often assumed that only a small number of predictors among
$X_1,\cdots,X_p$ contribute to the response, which amounts to
assuming ideally that the parameter vector $\bbeta$ is sparse. With
sparsity, variable selection can improve estimation accuracy by
effectively identifying the subset of important predictors, and also
enhance model interpretability with parsimonious representation.

Sparsity comes frequently with high dimensional data, which is a
growing feature in many areas of contemporary statistics. The
problems arise frequently in genomics such as gene expression and
proteomics studies, biomedical imaging, functional MRI, tomography,
tumor classifications, signal processing, image analysis, and
finance, where the number of variables or parameters $p$ can be much
larger than sample size $n$. For instance, one may wish to classify
tumors using microarray gene expression or proteomics data; one may
wish to associate protein concentrations with expression of genes or
predict certain clinical prognosis (e.g., injury scores or survival
time) using gene expression data. For this kind of problems, the
dimensionality can be much larger than the sample size, which calls
for new or extended statistical methodologies and theories. See,
e.g., Donoho (2000) and Fan and Li (2006) for overviews of
statistical challenges with high dimensionality.

Back to the problem in (\ref{011}), it is challenging to find tens
of important variables out of thousands of predictors, with number
of observations usually in tens or hundreds.  This is similar to
finding a couple of needles in a huge haystack. A new idea in Candes
and Tao (2007) is the notion of uniform uncertainty principle (UUP)
on deterministic design matrices. They proposed the Dantzig
selector, which is the solution to an $\ell_1$-regularization problem,
and showed that under UUP, this minimum $\ell_1$ estimator achieves
the ideal risk, i.e., the risk of the oracle estimator with the true
model known ahead of time, up to a logarithmic factor $\log p$.
Appealing features of the Dantzig selector include: 1) it is easy to
implement because the convex optimization the Dantzig selector
solves can easily be recast as a linear program; and 2) it has the
oracle property in the sense of Donoho and Johnstone (1994).

Despite their remarkable achievement, we still have four concerns
when the Dantzig selector is applied to high or ultra-high
dimensional problems. First, a potential hurdle is the computational
cost for large or huge scale problems such as implementing linear
programs in dimension tens or hundreds of thousands. Second, the
factor $\log p$ can become large and may not be negligible when
dimension $p$ grows rapidly with sample size $n$. Third, as
dimensionality grows, their UUP condition may be hard to satisfy,
which will be illustrated later using a simulated example. Finally,
there is no guarantee the Dantzig selector picks up the right model
though it has the oracle property. These four concerns inspire our
work.

\subsection{Dimensionality reduction}
Dimension reduction or feature selection is an effective strategy to
deal with high dimensionality. With dimensionality reduced from high to
low, computational burden can be reduced drastically. Meanwhile,
accurate estimation can be obtained by using some well-developed
lower dimensional method. Motivated by this along with those
concerns on the Dantzig selector, we have the following main goal in
our paper:
\begin{itemize}
\item Reduce dimensionality $p$ from a large or huge scale (say, $\exp(O(n^\xi))$
for some $\xi > 0$) to a relatively large scale $d$ (e.g., $o(n)$)
by a fast and efficient method.
\end{itemize}
We achieve this by introducing the concept of sure screening and
proposing a sure screening method based on a correlation learning
which filters out the features that have weak correlation with the
response.  Such a correlation screening is called Sure Independence
Screening (SIS). Here and below, by sure screening we mean a
property that all the important variables survive after variable
screening with probability tending to one. This dramatically narrows
down the search for important predictors. In particular, applying
the Dantzig selector to the much smaller submodel relaxes our first
concern on the computational cost. In fact, this not only speeds up
the Dantzig selector, but also reduces the logarithmic factor in
mimicking the ideal risk from $\log p$ to $\log d$, which is smaller
than $\log n$ and hence relaxes our second concern above.  It also
addresses he third concern since the UUP condition is easier to
satisfy.

Oracle properties in a stronger sense, say, mimicking the oracle in
not only selecting the right model, but also estimating the
parameters efficiently, give a positive answer to our third  and
fourth concerns above. Theories on oracle properties in this sense
have been developed in the literature. Fan and Li (2001) lay down
groundwork on variable selection problems in the finite parameter
setting. They discussed a family of variable selection methods that
adopt a penalized likelihood approach, which includes
well-established methods such as the AIC and BIC, as well as more
recent methods like the bridge regression in Frank and Friedman
(1993), Lasso in Tibshirani (1996), and SCAD in Fan (1997) and
Antoniadis and Fan (2001), and established oracle properties for
nonconcave penalized likelihood estimators. Later on, Fan and Peng
(2004) extend the results to the setting of $p=o(n^{1/3})$ and show
that the oracle properties continue to hold. An effective algorithm
for optimizing penalized likelihood, local quadratic approximation
(LQA), was proposed in Fan and Li (2001) and well studied in Hunter
and Li (2005). Zou (2006) introduces an adaptive Lasso in a finite
parameter setting and shows that Lasso does not have oracle
properties as conjectured in Fan and Li (2001), whereas the adaptive
Lasso does. Zou and Li (2008) propose a local linear approximation
algorithm that recasts the computation of non-concave penalized
likelihood problems into a sequence of penalized $L_1$-likelihood
problems.  They also proposed and studied the one-step sparse
estimators for nonconcave penalized likelihood methods.

There is a huge literature on the problem of variable selection. To
name a few in addition to those mentioned above, Fan and Li (2002)
study variable selection for Cox's proportional hazards model and
frailty model; Efron, Hastie, Johnstone and Tibshirani (2004)
propose LARS; Hunter and Li (2005) propose a new class of
algorithms, MM algorithms, for variable selection; Meinshausen and
B\"{u}hlmann (2006) look at the problem of variable selection with
the Lasso for high dimensional graphs, and Zhao and Yu (2006) give
an almost necessary and sufficient condition on model selection
consistency of Lasso. Meier, van de Geer and B\"uhlmann (2008)
 proposed a fast implementation for group Lasso. More recent studies
include Huang, Horowitz and Ma (2008), Paul {\em et al.} (2007),
Zhang (2007), and Zhang and Huang (2008), which signficantly
advances the theory and methods of the penalized least-squares
approaches. It is worth to mention that in variable selection, there
is a weaker concept than consistency, called persistency, introduced
by Greenshtein and Ritov (2004). Motivation of this concept lies in
the fact that in machine learning such as tumor classifications, the
primary interest centers on the misclassification errors or more
generally expected losses, not the accuracy of estimated parameters.
Greenshtein and Ritov (2004) study the persistency of Lasso-type
procedures in high dimensional linear predictor selection, and
Greenshtein (2006) extends the results to more general loss
functions. Meinshausen (2007) considers a case with finite
nonsparsity and shows that under quadratic loss, Lasso is
persistent, but the rate of persistency is slower than that of a
relaxed Lasso.

\subsection{Some insight on high dimensionality}
To gain some insight on challenges of high dimensionality in
variable selection, let us look at a situation where all the
predictors $X_1,\cdots,X_p$ are standardized and the distribution of
$\bz=\Sig^{-1/2}\bx$ is spherically symmetric, where
$\bx=\left(X_1,\cdots,X_p\right)\t$ and $\Sig=\cov\left(\bx\right)$.
Clearly, the transformed predictor vector $\bz$ has covariance
matrix $I_p$. Our way of study in this paper is to separate the
impacts of the covariance matrix $\Sig$ and the distribution of
$\bz$, which gives us a better understanding on difficulties of high
dimensionality in variable selection.

The real difficulty when dimension $p$ is larger than sample size
$n$ comes from four facts. First, the design matrix $\bX$ is
rectangular, having more columns than rows. In this case, the matrix
$\bX^T\bX$ is huge and singular. The maximum spurious correlation
between a covariate and the response can be large (see, e.g., Figure
1) because of the dimensionality and the fact that an unimportant
predictor can be highly correlated with the response variable due to
the presence of important predictors associated with the predictor.
These make variable selection difficult. Second, the population
covariance matrix $\Sig$ may become ill-conditioned as $n$ grows,
which adds difficulty to variable selection. Third, the minimum
nonzero absolute coefficient $|\beta_i|$ may decay with $n$ and get
close to the noise level, say, the order $(\log p / n)^{-1/2}$.
Fourth, the distribution of $\bz$ may have heavy tails. Therefore,
in general, it is challenging to estimate the sparse parameter
vector $\bbeta$ accurately when $p\gg n$.

\begin{figure} \centering
\begin{center}%
\includegraphics[trim=0.000000in 0.000000in 0.000000in -0.189256in,
height=2.25in, width=4in]%
{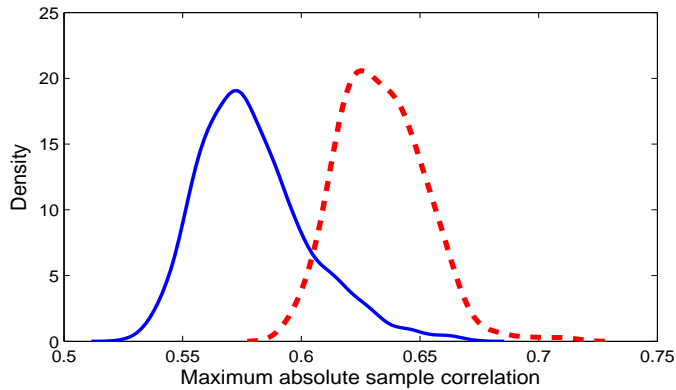}%
\\
\caption{Distributions of the maximum absolute sample correlation
coefficient
when $n=60, p=1000$ (solid curve) and $n=60, p=5000$ (dashed curve),
based on 500 simulations.}%
\end{center}%
\end{figure}%

When dimension $p$ is large, some of the intuition might not be
accurate. This is exemplified by the data piling problems in high
dimensional space observed in Hall, Marron and Neeman (2005). A
challenge with high dimensionality is that important predictors can
be highly correlated with some unimportant ones, which usually
increases with dimensionality. The maximum spurious correlation also
grows with dimensionality. We illustrate this using a simple
example. Suppose the predictors $X_1, \cdots, X_p$ are independent
and follow the standard normal distribution. Then, the design matrix
is an $n \times p$ random matrix, each entry an independent
realization from $\mathcal{N}(0,1)$. The maximum absolute sample
correlation coefficient among predictors can be very large. This is
indeed against our intuition, as the predictors are independent. To
show this, we simulated 500 data sets with $n=60$ and $p = 1000$ and
$p=5000$, respectively. Figure 1 shows the distributions of the
maximum absolute sample correlation. The multiple canonical
correlation between two groups of predictors (e.g., 2 in one group
and 3 in another) can even be much larger, as there are already ${p
\choose 2} {p-2 \choose 3}=O(p^5)$ choices of the two groups in our
example. Hence, sure screening when $p$ is large is very
challenging.

The paper is organized as follows. In the next section we propose a
sure screening method Sure Independence Screening (SIS) and discuss
its rationale as well as its connection with other methods of
dimensionality reduction. In Section 3 we review several known techniques
for model selection in the reduced feature space and present two
simulations and one real data example to study the performance of
SIS based model selection methods. In Section 4 we discuss some
extensions of SIS and in particular, an iterative SIS is proposed
and illustrated by three simulated examples. Section 5 is devoted to
the asymptotic analysis of SIS, an iteratively thresholded ridge
regression screener as well as two
SIS based model selection methods. Some concluding remarks are given
in Section 6. Technical details are provided in the Appendix.

\section{Sure Independence Screening}

\subsection{A sure screening method: correlation learning}
By sure screening we mean a property that all the important
variables survive after applying a variable screening procedure with
probability tending to one. A dimensionality reduction method is
desirable if it has the sure screening property. Below we introduce
a simple sure screening method using componentwise regression or
equivalently a correlation learning. Throughout the paper we center
each input variable so that the observed mean is zero, and scale
each predictor so that the sample standard deviation is one. Let
$\mathcal{M}_{*}=\left\{1\leq i\leq p :\beta_{i}\neq 0\right\}$ be
the true sparse model with nonsparsity size $s =|\mathcal{M}_{*}|$.
The other $p -s $ variables can also be correlated with the response
variable via linkage to the predictors contained in the model. Let
$\bomega=\left(\omega_1,\cdots,\omega_{p }\right)\t$ be a $p
$-vector obtained by the componentwise regression, that is,
\begin{equation} \label{120}
\bomega=\bX\t\by,
\end{equation}
where the $n\times p$ data matrix $\bX$ is first standardized
columnwise as mentioned before.  Hence, $\bomega$ is really a vector
of marginal correlations of predictors with the response variable,
rescaled by the standard deviation of the response.

For any given $\gamma \in\left(0,1\right)$, we sort the $p$
componentwise magnitudes of the vector $\bomega$ in a decreasing
order and define a submodel
\begin{equation} \label{110}
\mathcal{M}_{\gamma }=\left\{1\leq i\leq p
:\left|\omega_i\right|\text{ is among the first $\left[\gamma
n\right]$ largest of all}\right\},
\end{equation}
where $\left[\gamma n\right]$ denotes the integer part of $\gamma
n$. This is a straightforward way to shrink the full model
$\left\{1,\cdots,p \right\}$ down to a submodel $\mathcal{M}_{\gamma
}$ with size $d =[\gamma n]<n$. Such a correlation learning ranks
the importance of features according to their marginal correlation
with the response variable and filters out those that have weak
marginal correlations with the response variable.  We call this
correlation screening method Sure Independence Screening (SIS),
since each feature is used independently as a predictor to decide
how useful it is for predicting the response variable. This concept
is broader than the correlation screening and is applicable to
generalized linear models, classification problems under various
loss functions, and nonparametric learning under sparse additive
models.

The computational cost of correlation learning or SIS is that of
multiplying a $p\times n$ matrix with an $n$-vector plus getting the
largest $d$ components of a $p$-vector, so SIS has computational
complexity $O(np)$.

It is worth to mention that SIS uses only the order of componentwise
magnitudes of $\bomega$, so it is indeed invariant under scaling.
Thus the idea of SIS is identical to selecting predictors using
their correlations with the response. To implement SIS, we note that
linear models with more than $n$ parameters are not identifiable
with only $n$ data points.  Hence, we may choose $d=[\gamma n]$ to
be conservative, for instance, $n-1$ or $n /\log n$ depending on the
order of sample size $n$. Although SIS is proposed to reduce
dimensionality $p$ from high to below sample size $n$, nothing can
stop us applying it with final model size $d\geq n$, say,
$\gamma\geq1$. It is obvious that larger $d$ means larger
probability to include the true model $\mathcal{M}_*$ in the final
model $\mathcal{M}_\gamma$.

SIS is a hard-thresholding-type method. For orthogonal design
matrices, it is well understood. But for general design matrices,
there is no theoretical support for it, though this kind of idea is
frequently used in applications. It is important to identify the
conditions under which the sure screening property holds for SIS,
i.e.,
\begin{equation} \label{121}
P\left(\mathcal{M}_{*}\subset\mathcal{M}_{\gamma
}\right)\rightarrow1\quad\text{as }n\rightarrow\infty
\end{equation}
for some given $\gamma $. This question as well as how the sequence
$\gamma=\gamma_n \to 0$ should be chosen will be answered by Theorem
1 in Section 5. We would like to point out that the Simple
Thresholding Algorithm (see, e.g., Baron \textit{et al.}, 2005 and
Gribonval \textit{et al.}, 2007) that is used in sparse
approximation or compressed sensing is a one step greedy algorithm
and related to SIS. In particular, our asymptotic analysis in
Section 5 helps to understand the performance of the Simple
Thresholding Algorithm.

\subsection{Rationale of correlation learning}
To better understand the rationale of the correlation learning, we
now introduce an iteratively thresholded ridge regression screener
(ITRRS), which is an extension of the dimensionality reduction
method SIS. But for practical implementation, only the correlation
learning is needed. ITRRS also provides a very nice technical tool
for our understanding of the sure screening property of the
correlation screening and other methods.

When there are more predictors than observations, it is well known
that the least squares estimator
$\hbbeta_{\text{LS}}=\left(\bX\t\bX\right)^{+}\bX\t\by$ is noisy,
where $\left(\bX\t\bX\right)^{+}$ denotes the Moore-Penrose
generalized inverse of $\bX\t\bX$. We therefore consider the ridge
regression, namely, linear regression with $\ell_2$-regularization
to reduce the variance. Let $\bomega^{\lambda }=(\omega_1^{\lambda
},\cdots,\omega_{p }^{\lambda })\t$ be a $p $-vector obtained by the
ridge regression, that is,
\begin{equation} \label{012}
\bomega^{\lambda }=\left(\bX\t\bX+\lambda I_{p
}\right)^{-1}\bX\t\by,
\end{equation}
where $\lambda >0$ is a regularization parameter. It is obvious that
\begin{equation} \label{118}
\bomega^{\lambda }\rightarrow\hbbeta_{\text{LS}}\quad \text{as }
\lambda \rightarrow0,
\end{equation}
and the scaled ridge regression estimator tends to the componentwise
regression estimator:
\begin{equation} \label{096}
\lambda \bomega^{\lambda }\rightarrow\bomega\quad \text{as } \lambda
\rightarrow\infty.
\end{equation}
In view of (\ref{118}), to make $\bomega^{\lambda }$ less noisy we
should choose large regularization parameter $\lambda $ to reduce
the variance in the estimation.  Note that the ranking of the
absolute components of $\bomega^{\lambda }$ is the same as that of
$\lambda \bomega^{\lambda }$.  In light of (\ref{096}) the
componentwise regression estimator is a specific case of the ridge regression
with regularization parameter $\lambda = \infty$, namely, it makes the resulting
estimator as less noisy as possible.

For any given $\delta \in\left(0,1\right)$, we sort the $p$
componentwise magnitudes of the vector $\bomega^{\lambda }$ in a
descending order and define a submodel
\begin{equation} \label{039}
\mathcal{M}^{1}_{\delta, \lambda}=\left\{1\leq i\leq p
:|\omega_i^{\lambda }|\text{ is among the first $\left[\delta p
\right]$ largest of all}\right\}.
\end{equation}
This procedure reduces the model size by a factor of $1-\delta $.
The idea of ITRRS to be introduced below is to perform
dimensionality reduction as above successively until the number of
remaining variables drops to below sample size $n$.

It will be shown in Theorem 2 in Section 5 that under some
regularity conditions and when the tuning parameters $\lambda$ and
$\delta$ are chosen appropriately, with overwhelming probability the
submodel $\mathcal{M}^{1}_{\delta, \lambda}$ will contain the true
model $\mathcal{M}_{*}$ and its size is an order $n^\theta$ for some
$\theta > 0$ lower than the original one $p $. This property
stimulates us to propose ITRRS as follows:
\begin{itemize}
\item First, carry out the
procedure in (\ref{039}) to the full model $\left\{1,\cdots,p
\right\}$ and get a submodel $\mathcal{M}^{1}_{\delta, \lambda}$
with size $[\delta p ]$;

\item Then, apply a similar procedure to the model
$\mathcal{M}^{1}_{\delta, \lambda}$ and again obtain a submodel
$\mathcal{M}^{2}_{\delta, \lambda}\subset\mathcal{M}^{1}_{\delta, \lambda}$
with size $[\delta ^2p ]$, and so on;

\item Finally, get a submodel
$\mathcal{M}_{\delta, \lambda}=\mathcal{M}^{k}_{\delta, \lambda}$
with size $d =[\delta ^{k }p ]<n$, where $[\delta ^{k -1}p ]\geq n$.
\end{itemize}
We would like to point out that the above procedure is different
from the threshholded ridge regression, as the submodels and
estimated parameters change over the course of iterations.  The only
exception is the case that $\lambda = \infty$, in which the rank of
variables do not vary with iterations.

Now we are ready to see that the correlation learning introduced in
Section 2.1 is a specific case of ITRRS since the componentwise
regression is a specific case of the ridge regression with an
infinite regularization parameter. The ITRRS provides a very nice
technical tool for understanding how fast the dimension $p$ can grow
compared with sample size $n$ and how the final model size $d$ can
be chosen while the sure screening property still holds for the
correlation learning. The question of whether ITRRS has the sure
screening property as well as how the tuning parameters $\gamma$ and
$\delta$ should be chosen will be answered by Theorem 3 in Section
5.

The number of steps in ITRRS depends on the choice of
$\delta\in(0,1)$. We will see in Theorem 3 that $\delta$ can not be
chosen too small which means that there should not be too many
iteration steps in ITRRS. This is due to the cumulation of the
probability errors of missing some important variables over the
iterations. In particular, the backward stepwise deletion regression
which deletes one variable each time in ITRRS until the number of
remaining variables drops to below sample size might not work in
general as it requires $p-d$ iterations. When $p$ is of exponential
order, even though the probability of mistakenly deleting some
important predictors in each step of deletion is exponentially
small, the cumulative error in exponential order of operations may
not be negligible.

\subsection{Connections with other dimensionality reduction methods}
As pointed out before, SIS uses the marginal information of
correlation to perform dimensionality reduction. The idea of using
marginal information to deal with high dimensionality has also
appeared independently in Huang, Horowitz and Ma (2008) who proposed
to use marginal bridge estimators to select variables for sparse
high dimensional regression models. We now look at SIS in the
context of classification, in which the idea of independent
screening appears natural and has been widely used.

The problem of classification can be regarded as a specific case of
the regression problem with response variable taking discrete values
such as $\pm1$. For high dimensional problems like tumor
classification using gene expression or proteomics data, it is not
wise to classify the data using the full feature space due to the
noise accumulation and interpretability.  This is well demonstrated
both theoretically and numerically in Fan and Fan (2008). In
addition, many of the features come into play through linkage to the
important ones (see, e.g., Figure 1). Therefore feature selection is
important for high dimensional classification. How to effectively
select important features and how many of them to include are two
tricky questions to answer. Various feature selection procedures
have been proposed in the literature to improve the classification
power in presence of high dimensionality. For example, Tibshirani
\textit{et al.} (2002) introduce the nearest shrunken centroids
method, and Fan and Fan (2008) propose the Features Annealed
Independence Rules (FAIR) procedure. Theoretical justification for
these methods are given in Fan and Fan (2008).

SIS can readily be used to reduce the feature space. Now suppose we
have $n_1$ samples from class $1$ and $n_2$ samples from class $-1$.
Then the componentwise regression estimator (\ref{120}) becomes
\begin{equation}
\bomega = \sum_{Y_i = 1} \bx_i  - \sum_{Y_i = -1} \bx_i.,
\label{a01}
\end{equation}
Written more explicitly, the $j$-th component of the $p$-vector
$\bomega$ is
$$
\omega_j = (n_1\bar{X}_{j,1} - n_2 \bar{X}_{j, 2})/\mbox{SD of the
$j$-th feature},
$$
by recalling that each covariate in (\ref{a01}) has been normalized
marginally, where $\bar{X}_{j,1}$ is the sample average of the
$j$-th feature with class label ``1'' and $\bar{X}_{j,2}$ is the
sample average of the $j$-th feature with class label ``$-1$''. When
$n_1 = n_2$, $\omega_j$ is simply a version of the two-sample
$t$-statistic except for a scaling constant. In this case, feature
selection using SIS is the same as that using the two-sample
$t$-statistics.  See Fan and Fan (2008) for a theoretical study of
sure screening property in this context.

Two-sample $t$-statistics are commonly used in feature selection for
high dimensional classification problems such as in the significance
analysis of gene selection in microarray data analysis (see, e.g.,
Storey and Tibshirani, 2003; Fan and Ren, 2006) as well as in the
nearest shrunken centroids method of Tibshirani \textit{et al.}
(2002). Therefore SIS is an insightful and natural extension of this
widely used technique. Although not directly applicable, the sure
screening property of SIS in Theorem 1 after some adaptation gives
theoretical justification for the nearest shrunken centroids method.
See Fan and Fan (2008) for a sure screening property.

By using SIS we can single out the important features and thus
reduce significantly the feature space to a much lower dimensional
one. From this point on, many methods such as the linear
discrimination (LD) rule or the naive Bayes (NB) rule can be applied
to conduct the classification in the reduced feature space. This
idea will be illustrated on a Leukemia data set in Section 3.3.3.

\section{SIS based model selection techniques}
\subsection{Estimation and model selection in the reduced feature space}
As shown later in Theorem 1 in Section 5, with the correlation
learning, we can shrink the full model $\left\{1,\cdots,p \right\}$
straightforward and accurately down to a submodel
$\mathcal{M}=\mathcal{M}_{\gamma }$ with size $d =[\gamma n]=o(n)$.
Thus the original problem of estimating the sparse $p $-vector
$\bbeta $ in (\ref{011}) reduces to estimating a sparse $d $-vector
$\bbeta =\left(\beta_{1},\cdots,\beta_{d }\right)\t$ based on the
now much smaller submodel $\mathcal{M}$, namely,
\begin{equation} \label{001}
\by=\bX_{\mathcal{M}}\bbeta +\bveps,
\end{equation}
where $\bX_{\mathcal{M}}=\left(\bx_1,\cdots,\bx_n\right)\t$ denotes
an $n\times d $ submatrix of $\bX$ obtained by extracting its
columns corresponding to the indices in $\mathcal{M}$. Apparently
SIS can speed up variable selection dramatically when the original
dimension $p$ is ultra high.

Now we briefly review several well-developed moderate dimensional
techniques that can be applied to estimate the $d $-vector $\bbeta $
in (\ref{001}) at the scale of $d$ that is comparable with $n$.
Those methods include SCAD in Fan and Li (2001) and Fan and Peng
(2004), adaptive Lasso in Zou (2006), the Dantzig selector in Candes
and Tao (2007), among others.

\subsubsection{Penalized least-squares and SCAD}
Penalization is commonly used in variable selection. Fan and Li
(2001, 2006) give a comprehensive overview of feature selection and
a unified framework based on penalized likelihood approach to the
problem of variable selection. They consider the penalized least
squares (PLS)
\begin{equation} \label{111}
\ell\left(\bbeta\right)=\frac{1}{2n}\sum_{i=1}^n\left(Y_i-\bx_i\t\bbeta\right)^2+
\sum_{j=1}^{d }p_{\lambda_j}(\left|\beta_j\right|),
\end{equation}
where $\bbeta=\left(\beta_1,\cdots,\beta_{d
}\right)\t\in\mathbf{R}^{d }$ and $p_{\lambda_j}(\cdot)$ is a
penalty function indexed by a regularization parameter $\lambda_j$.
Variation of the regularization parameters across the predictors
allows us to incorporate some prior information. For example, we may
want to keep certain important predictors in the model and choose
not to penalize their coefficients. The regularization parameters
$\lambda_j$ can be chosen, for instance, by cross-validation (see,
e.g., Breiman, 1996 and Tibshirani, 1996). A unified and effective
algorithm for optimizing penalized likelihood, called local
quadratic approximation (LQA), was proposed in Fan and Li (2001) and
well studied in Hunter and Li (2005). In particular, LQA can be
employed to minimize the above PLS.  In our implementation,
we choose $\lambda_j = \lambda$ and select $\lambda$ by BIC.

\begin{figure} \centering
\begin{center}%
\includegraphics[height=2.25in] {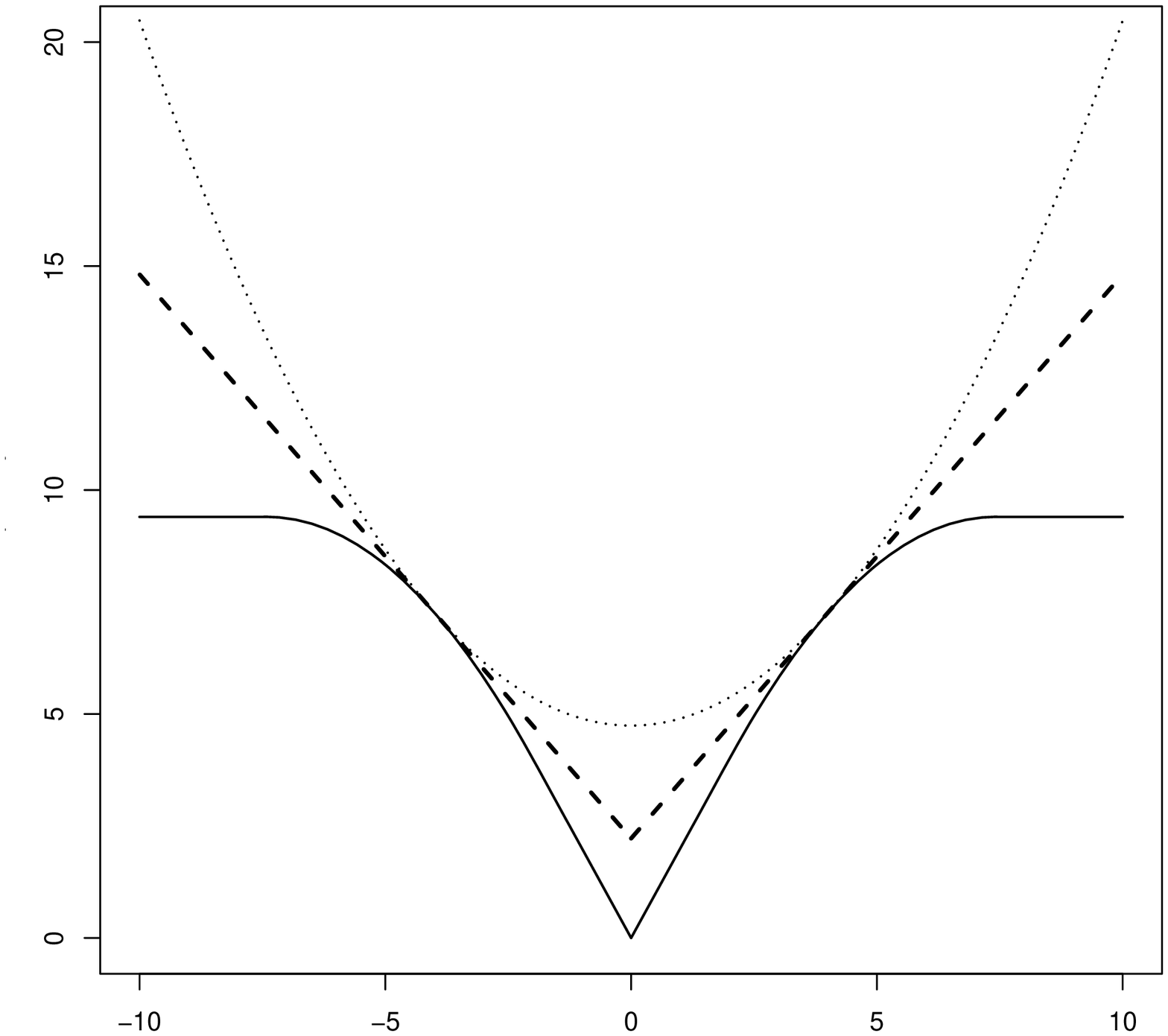} \hspace*{0.5 in}
\includegraphics[height=2.35in] {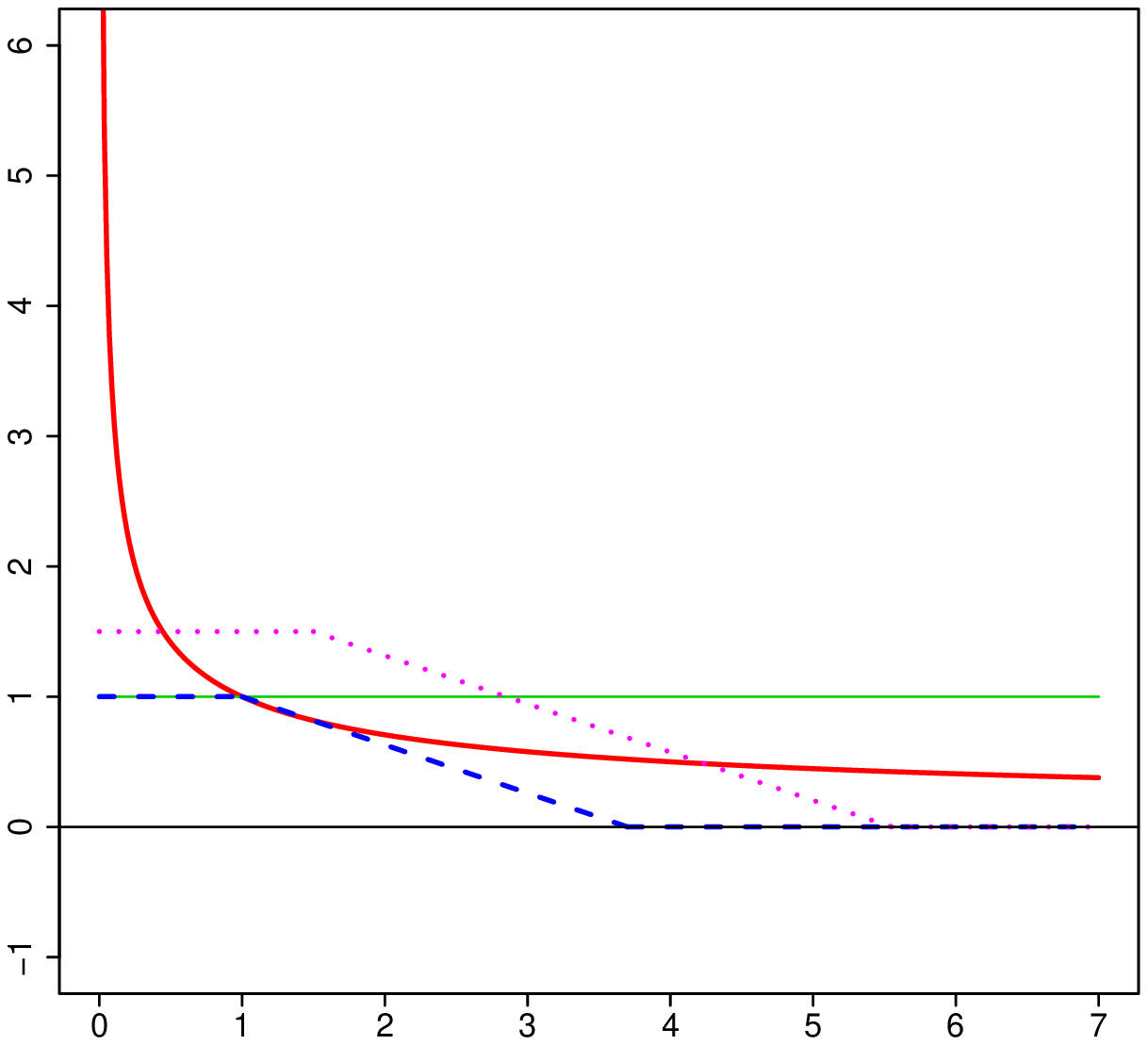}
\\
\caption{Left panel:  The SCAD penalty (solid) and its local linear (dashed)
 and quadratic (dotted) approximations at the point $x=4$.  Right panel:
 $p'_\lambda(\cdot)$ for penalized $L_1$ (thin solid), SCAD with $\lambda = 1$ (dashed)
 and $\lambda = 1.5$ (dotted) and adaptive Lasso (thick solid) with $\gamma = 0.5$.}%
\end{center}%
\end{figure}%

An alternative and effective algorithm to minimize the penalized
least-squares problem (\ref{111}) is the local linear approximation
(LLA) proposed by Zou and Li (2008).  With the local linear
approximation, the problem (\ref{111}) can be cast as a sequence of
penalized $L_1$ regression problems so that the LARS (Efron, {\em et
al.}, 2004) or other algorithms can be employed. More explicitly,
given the estimate $\{\hat{\beta}_{j}^{(k)}, j = 1, \cdots, d \}$ at
the $k$-th iteration, instead of minimizing (\ref{111}), one
minimizes
\begin{equation}
   \label{jf1} \frac{1}{2n}\sum_{i=1}^n\left(Y_i-\bx_i\t\bbeta\right)^2+
\sum_{j=1}^{d }  w_j^{(k)} |\beta_j|,
\end{equation}
which after adding the constant term $\sum_{j=1}^{d
}p_{\lambda_j}(|\hat{\beta}_j^{(k)} |)$ is a local linear
approximation to $\ell(\bbeta)$ in (\ref{111}), where $w_j^{(k)} =
|p_{\lambda_j}'(|\hat{\beta}_j^{(k)} |) |$. Problem (\ref{jf1}) is a
convex problem and can be solved by LARS and other algorithm such as
those in Friedman {\em et al.} (2007) and Meier, van der Geer and
B\"uhlmann (2008). In this sense, the penalized least-squares
problem (\ref{111}) can be regarded as a family of weighted
penalized $L_1$-problem and the function $p_\lambda'(\cdot)$
dictates the amount of penalty at each location.  The emphasis on
non-concave penalty functions by Fan and Li (2001) is to ensure that
penalty decreases to zero as $|\hat{\beta}_j^{(k)} |$ gets large.
This reduces unnecessary biases of the penalized likelihood
estimator, leading to the oracle property in Fan and Li (2001).
Figure 2 depicts how the SCAD function is approximated locally by a
linear or quadratic function and the derivative functions
$p'_\lambda(\cdot)$ for some commonly used penalty functions. When
the initial value $\bbeta = 0$, the first step estimator is indeed
LASSO so the implementation of SCAD can be regarded as an
iteratively reweighted penalized $L_1$-estimator with LASSO as an
initial estimator. See Section 6 for further discussion of the
choice of initial values $\{\hat{\beta}_{j}^{(0)}, j = 1, \cdots, d
\}$.

The PLS (\ref{111}) depends on the choice of penalty function
$p_{\lambda_j}(\cdot)$. Commonly used penalty functions include the
$\ell_p$-penalty, $0\leq p\leq 2$, nonnegative garrote in Breiman
(1995), and smoothly clipped absolute deviation (SCAD) penalty, in
Fan (1997) and a minimax concave penality (MCP) in Zhang (2007) (see
below for definition). In particular, the $\ell_1$-penalized least
squares is called Lasso in Tibshirani (1996). In seminal papers,
Donoho and Huo (2001) and Donoho and Elad (2003) show that penalized
$\ell_0$-solution can be found by penalized $\ell_1$-method when the
problem is sparse enough, which implies that the best subset
regression can be found by using the penalized $\ell_1$-regression.
Antoniadis and Fan (2001) propose the PLS for wavelets denoising
with irregular designs. Fan and Li (2001) advocate penalty functions
with three properties: sparsity, unbiasedness, and continuity. More
details on characterization of these three properties can be found
in Fan and Li (2001) and Antoniadis and Fan (2001). For penalty
functions, they showed that singularity at the origin is a necessary
condition to generate sparsity and nonconvexity is required to
reduce the estimation bias. It is well known that $\ell_p$-penalty
with $0\leq p<1$ does not satisfy the continuity condition,
$\ell_p$-penalty with $p>1$ does not satisfy the sparsity condition,
and $\ell_1$-penalty (Lasso) possesses the sparsity and continuity,
but generates estimation bias, as demonstrated in Fan and Li (2001),
Zou (2006), and Meinshausen (2007).

Fan (1997) proposes a continuously differentiable penalty function
called the smoothly clipped absolute deviation (SCAD) penalty, which
is defined by
\begin{equation} \label{112}
p'_\lambda(\left|\beta\right|)=\lambda\left\{I\left(\left|\beta\right|\leq\lambda\right)
+\frac{\left(a\lambda-\left|\beta\right|\right)_+}{\left(a-1\right)\lambda}
I\left(\left|\beta\right|>\lambda\right)\right\}\quad\text{for some
}a>2.
\end{equation}
Fan and Li (2001) suggest using $a=3.7$. This function has similar
feature to the penalty function
$\lambda\left|\beta\right|/\left(1+\left|\beta\right|\right)$
advocated in Nikolova (2000). The MCP in Zhang (2007) translates the
flat part of the derivative of the SCAD to the origin and is given
by
\[ p_\lambda'(|\beta|) = (a \lambda - |\beta|)_+/a, \]
which minimizes the maximum of the concavity. The SCAD penalty and
MCP satisfy the above three conditions simultaneously. We will show
in Theorem 5 in Section 5 that SIS followed by the SCAD enjoys the
oracle properties.

\subsubsection{Adaptive Lasso}
The Lasso in Tibshirani (1996) has been widely used due to its
convexity. It however generates estimation bias. This problem was
pointed out in Fan and Li (2001) and formally shown in Zou (2006)
even in a finite parameter setting. To overcome this bias problem,
Zou (2006) proposes an adaptive Lasso and Meinshausen (2007)
proposes a relaxed Lasso.

The idea in Zou (2006) is to use an adaptively weighted $\ell_1$
penalty in the PLS (\ref{111}). Specifically, he introduced the
following penalization term
\[ \lambda \sum_{j=1}^d\omega_j\left|\beta_j\right|, \]
where $\lambda \geq0$ is a regularization parameter and
$\bomega=\left(\omega_1,\cdots,\omega_d\right)\t$ is a known weight
vector. He further suggested using the weight vector
$\widehat{\bomega}=1/|\hbbeta|^\gamma$, where $\gamma\geq0$, the
power is understood componentwise, and $\hbbeta$ is a root-$n$
consistent estimator. In view of (\ref{jf1}), the adaptive Lasso is
really the implementation of PLS (\ref{111}) with
$p_\lambda(|\beta|) = |\beta|^{1-\gamma}$ using LLA. Its connections
with the family of non-concave penalized least-squares is apparently
from (\ref{jf1}) and Figure 2.

The case of $\gamma=1$ is closely related to the nonnegative garrote
in Breiman (1995). Zou (2006) also showed that the adaptive Lasso
can be solved by the LARS algorithm, which was proposed in Efron,
Hastie, Johnstone and Tibshirani (2004). Using the same finite
parameter setup as that in Knight and Fu (2000), Zou (2006)
establishes that the adaptive Lasso has the oracle properties as
long as the tuning parameter is chosen in a way such that $\lambda
/\sqrt{n}\rightarrow0$ and $\lambda
n^{\frac{\gamma-1}{2}}\rightarrow\infty$ as $n\rightarrow\infty$.

\subsubsection{Dantzig selector}
The Dantzig selector was proposed in Candes and Tao (2007) to
recover a sparse high dimensional parameter vector in the linear
model. Adapted to the setting in (\ref{001}), it is the solution
$\hbbeta_{\text{DS}}$ to the following $\ell_1$-regularization
problem
\begin{equation} \label{002}
\min_{\bzeta\in\mathbf{R}^{d
}}\left\|\bzeta\right\|_1\quad\text{subject to
}\left\|(\bX_\mathcal{M})\t\br\right\|_\infty\leq\lambda_{d }\sigma,
\end{equation}
where $\lambda_{d }>0$ is a tuning parameter,
$\br=\by-\bX_\mathcal{M}\bzeta$ is an $n$-vector of the residuals,
and $\|\cdot\|_1$ and $\|\cdot\|_\infty$ denote the $\ell_1$ and
$\ell_\infty$ norms, respectively. They pointed out that the above
convex optimization problem can easily be recast as a linear
program:
\[ \min\sum_{i=1}^{d }u_i\quad\text{subject
to }-\bu\leq\bzeta\leq\bu\text{ and }-\lambda_{d }\sigma\bone\leq
(\bX_\mathcal{M})\t\left(\by-\bX_\mathcal{M}\bzeta\right)\leq\lambda_{d
}\sigma\bone,
\] where the optimization variables are
$\bu=\left(u_1,\cdots,u_{d }\right)\t$ and $\bzeta\in\mathbf{R}^{d
}$, and $\bone$ is a $d $-vector of ones.

We will show in Theorem 4 in Section 5 that an application of SIS
followed by the Dantzig selector can achieve the ideal risk up to a
factor of $\log d $ with $d<n$, rather than the original $\log p$.
In particular, if dimension $p$ is growing exponentially fast, i.e.,
$p=\exp(O(n^\xi))$ for some $\xi > 0$, then a direct application of
the Dantzig selector results in a loss of a factor $O(n^\xi)$ which
could be too large to be acceptable. On the other hand, with the
dimensionality first reduced by SIS the loss is now merely of a
factor $\log d $, which is less than $\log n$.

\begin{figure} [htb]\centering
\begin{center}%
\includegraphics[trim=0.500000in 0.000000in 0.000000in -0.189256in,
height=1.35in, width=6.5in]%
{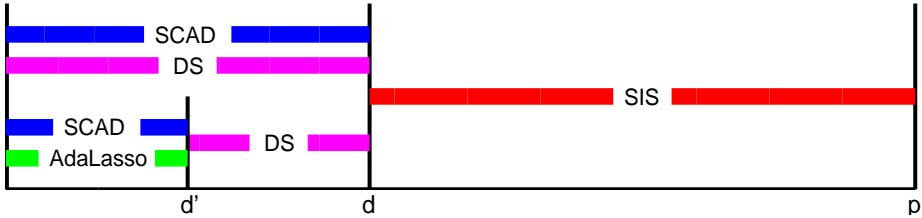}%
\\
\caption{Methods of model selection with ultra high dimensionality.}%
\end{center}%
\end{figure}%

\subsection{SIS based model selection methods}
For the problem of ultra-high dimensional variable selection, we
propose first to apply a sure screening method such as SIS to reduce
dimensionality from $p $ to a relatively large scale $d $, say,
below sample size $n$. Then we use a lower dimensional model
selection method such as the SCAD, Dantzig selector, Lasso, or
adaptive Lasso. We call SIS followed by the SCAD and Dantzig
selector SIS-SCAD and SIS-DS, respectively for short in the paper.
In some situations, we may want to further reduce the model size
down to $d'<d$ using a method such as the Dantzig selector along
with the hard thresholding or the Lasso with a suitable tuning, and
finally choose a model with a more refined method such as the SCAD
or adaptive Lasso. In the paper these two methods will be referred
to as SIS-DS-SCAD and SIS-DS-AdaLasso, respectively for simplicity.
Figure 3 shows a schematic diagram of these approaches.

The idea of SIS makes it feasible to do model selection with ultra
high dimensionality and speeds up variable selection drastically. It
also makes the model selection problem efficient and modular. SIS
can be used in conjunction with any model selection technique
including the Bayesian methods (see, e.g., George and McCulloch,
1997) and Lasso.  We did not include SIS-Lasso for numerical studies
due to the approximate equivalence between Dantzig selector and
Lasso (Bickel, Ritov and Tsybakov, 2007; Meinshausen, Rocha and Yu,
2007).

\subsection{Numerical studies}
To study the performance of SIS based model selection methods
proposed above, we now present two simulations and one real data
example.

\subsubsection{Simulation I: ``independent" features}
For the first simulation, we used the linear model (\ref{011}) with
i.i.d. standard Gaussian predictors and Gaussian noise with standard
deviation $\sigma=1.5$. We considered two such models with
$(n,p)=(200,1000)$ and $(800,20000)$, respectively. The sizes $s$ of
the true models, i.e., the numbers of nonzero coefficients, were
chosen to be 8 and 18, respectively, and the nonzero components of
the $p$-vectors $\bbeta $ were randomly chosen as follows. We set $a
=4\log n/\sqrt{n}$ and $5\log n/\sqrt{n}$, respectively, and picked
nonzero coefficients of the form $\left(-1\right)^u\left(a
+\left|z\right|\right)$ for each model, where $u$ was drawn from a
Bernoulli distribution with parameter $0.4$ and $z$ was drawn from
the standard Gaussian distribution. In particular, the
$\ell_2$-norms $\|\bbeta \|$ of the two simulated models are 6.795
and 8.908, respectively. For each model we simulated 200 data sets.
Even with i.i.d. standard Gaussian predictors, the above settings
are nontrivial since there is nonnegligible sample correlation among
the predictors, which reflects the difficulty of high dimensional
variable selection. As an evidence, we report in Figure 4 the
distributions of the maximum absolute sample correlation when
$n=200$ and $p=1000$ and $5000$, respectively. It reveals
significant sample correlation among the predictors. The multiple
canonical correlation between two groups of predictors can be much
larger.

\begin{figure} \centering
\begin{center}%
\includegraphics[trim=0.000000in 0.000000in 0.000000in -0.189256in,
height=2.25in, width=4in]%
{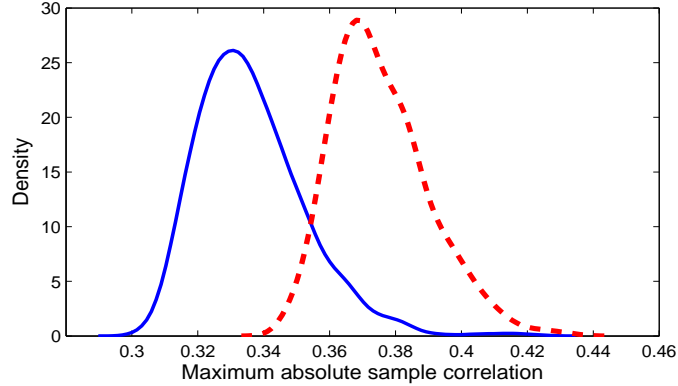}%
\\
\caption{Distributions of the maximum absolute sample correlation
when $n=200, p =1000$ (solid curve) and $n=200, p =5000$ (dashed curve).}%
\end{center}%
\end{figure}%

To estimate the sparse $p$-vectors $\bbeta$, we employed six
methods: the Dantzig selector (DS) using a primal-dual algorithm,
Lasso using the LARS algorithm, SIS-SCAD, SIS-DS, SIS-DS-SCAD, and
SIS-DS-AdaLasso (see Figure 3). For SIS-SCAD and SIS-DS, we chose $d
=[n/\log n]$ and for the last two methods, we chose $d=n-1$ and
$d'=[n/\log n]$ and in the middle step the Dantzig selector was used
to further reduce the model size from $d$ to $d'$ by choosing
variables with the $d'$ largest componentwise magnitudes of the
estimated $d$-vector (see Figure 3).

\begin{figure} \centering
\begin{center}%
\includegraphics[trim=0.000000in 0.000000in 0.000000in -0.189256in,
height=1.9in, width=6in]%
{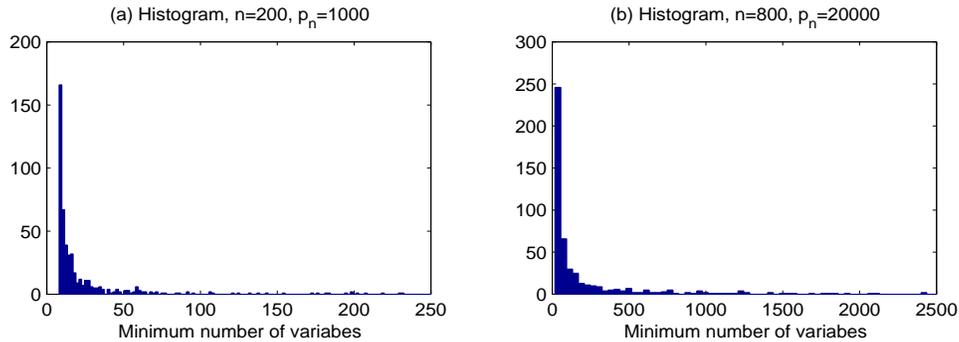}%
\\
\caption{{(a) Distribution of the minimum number of selected
variables required to include the true model by using SIS when
$n=200,p =1000$ in simulation I. (b) The same plot when $n=800,p=20000$.}}%
\end{center}%
\end{figure}%

\begin{table}[tph]
\begin{center}
Table 1: Results of simulation I\vspace{0.2 in}\\\small
\begin{tabular}
[c]{cccccccccccc}\hline  & \multicolumn{11}{c}{\bf Medians of the
selected model sizes (upper entry) }\\
& \multicolumn{11}{c}{\bf and the estimation errors (lower entry)}
\\\cline{2-12}
$p $ & \bf DS & & \bf Lasso & & \bf SIS-SCAD & & \bf SIS-DS & &
\bf SIS-DS-SCAD & & \bf SIS-DS-AdaLasso \\
\hline%
$1000$ & $10^3$ & & 62.5 & & 15 & &
37 & & 27 & & 34 \\
\cline{2-2} \cline{4-4} \cline{6-6} \cline{8-8} \cline{10-10}
\cline{12-12}
& 1.381 & & 0.895 & & 0.374 & & 0.795 & & 0.614 & & 1.269\\
20000 & --- & & --- & & 37 & &
119 & & 60.5 & & 99 \\
\cline{2-2} \cline{4-4} \cline{6-6} \cline{8-8} \cline{10-10}
\cline{12-12}
& --- & & --- & & 0.288 & & 0.732 & & 0.372 & & 1.014\\
\hline
\end{tabular}
\end{center}
\end{table}

The simulation results are summarized in Figure 5 and Table 1.
Figure 5, produced based on 500 simulations, depicts the
distribution of the minimum number of selected variables, i.e., the
selected model size, that is required to include all variables in
the true model by using SIS. It shows clearly that in both settings
it is safe to shrink the full model down to a submodel of size
$[n/\log n]$ with SIS, which is consistent with the sure screening
property of SIS shown in Theorem 1 in Section 5. For example, for
the case of $n = 200$ and $p = 1000$, reducing the model size to 50
includes the variables in the true model with high probability, and
for the case of $n = 800$ and $p  = 20000$, it is safe to reduce the
dimension to about 500. For each of the above six methods, we report
in Table 1 the median of the selected model sizes and median of the
estimation errors $\|\hbbeta -\bbeta \|$ in $\ell_2$-norm. Four
entries of Table 1 are missing due to limited computing power and
software used. In comparison, SIS reduces the computational burden
significantly.

From Table 1 we see that the Dantzig selector gives nonsparse
solutions and the Lasso using the cross-validation for selecting its
tuning parameter produces large models.   This can be due to the
fact that the biases in Lasso require a small bandwidth in
cross-validation, whereas a small bandwidth results in lack of
sparsistency, using the terminology of Ravikumar {\em et al.}
(2007). This has also been observed and demonstrated in the work by
Lam and Fan (2007) in the context of estimating sparse covariance or
precision matrices. We should point out here that a variation of the
Dantzig selector, the Gauss-Dantzig selector in Candes and Tao
(2007), should yield much smaller models, but for simplicity we did
not include it in our simulation. Among all methods, SIS-SCAD
performs the best and generates much smaller and more accurate
models. It is clear to see that SCAD gives more accurate estimates
than the adaptive Lasso in view of the estimation errors. Also, SIS
followed by the Dantzig selector improves the estimation accuracy
over using the Dantzig selector alone, which is in line with our
theoretical result.

\begin{figure} \centering
\begin{center}%
\includegraphics[trim=0.7000000in 0.000000in 0.000000in -0.189256in,
height=1.7in, width=6.6in]%
{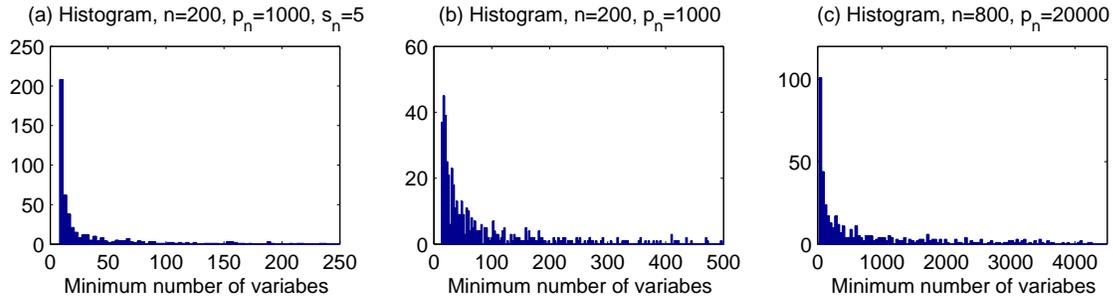}%
\\
\caption{{(a) Distribution of the minimum number of selected
variables required to include the true model by using SIS when
$n=200,p =1000,s =5$ in simulation II. (b) The same plot when
$n=200,p =1000,s =8$. (c) The same plot when $n=800,p =20000$.}}%
\end{center}%
\end{figure}%

\begin{table}[tph]
\begin{center}
Table 2: Results of simulation II\vspace{0.2 in}\\\small
\begin{tabular}
[c]{cccccccccccc}\hline  & \multicolumn{11}{c}{\bf Medians of the
selected model sizes (upper entry) }\\
& \multicolumn{11}{c}{\bf and the estimation errors (lower entry)}
\\\cline{2-12}
$p $ & \bf DS & & \bf Lasso & & \bf SIS-SCAD & & \bf SIS-DS & &
\bf SIS-DS-SCAD & & \bf SIS-DS-AdaLasso \\
\hline%
$1000$ & $10^3$ & & 91 & & 21 & &
56 & & 27 & & 52 \\
\cline{2-2} \cline{4-4} \cline{6-6} \cline{8-8} \cline{10-10}
\cline{12-12}
($s =5$)& 1.256 & & 1.257 & & 0.331 & & 0.727 & & 0.476 & & 1.204\\
 & $10^3$ & & 74 & & 18 & &
56 & & 31.5 & & 51 \\
\cline{2-2} \cline{4-4} \cline{6-6} \cline{8-8} \cline{10-10}
\cline{12-12}
($s =8$) & 1.465 & & 1.257 & & 0.458 & & 1.014 & & 0.787 & & 1.824\\
20000 & --- & & --- & & 36 & &
119 & & 54 & & 86 \\
\cline{2-2} \cline{4-4} \cline{6-6} \cline{8-8} \cline{10-10}
\cline{12-12}
& --- & & --- & & 0.367 & & 0.986 & & 0.743 & & 1.762\\
\hline
\end{tabular}
\end{center}
\end{table}

\subsubsection{Simulation II: ``dependent" features}
For the second simulation, we used similar models to those in
simulation I except that the predictors are now correlated with each
other. We considered three models with $(n,p,s)=(200, 1000, 5)$,
$(200, 1000, 8)$, and $(800, 20000, 14)$, respectively, where $s$
denotes the size of the true model, i.e., the number of nonzero
coefficients. The three $p$-vectors $\bbeta $ were generated in the
same way as in simulation I. We set $(\sigma,a)=(1,2\log
n/\sqrt{n})$, $(1.5,4\log n/\sqrt{n})$, and $(2,4\log n/\sqrt{n})$,
respectively. In particular, the $\ell_2$-norms $\|\bbeta \|$ of the
three simulated models are 3.304, 6.795, and 7.257, respectively. To
introduce correlation between predictors, we first used a Matlab
function \verb=sprandsym= to randomly generate an $s \times s $
symmetric positive definite matrix $\bA$ with condition number
$\sqrt{n}/\log n$, and drew samples of $s $ predictors
$X_1,\cdots,X_{s }$ from $\mathcal{N}(\bzero,\bA)$. Then we took
$Z_{s +1},\cdots,Z_{p }\sim\mathcal{N}(\bzero,I_{p -s })$ and
defined the remaining predictors as $X_i=Z_i+rX_{i-s }$, $i=s
+1,\cdots,2s $ and $X_i=Z_i+\left(1-r\right)X_1$, $i=2s +1,\cdots,p
$ with $r=1-4\log n/p $, $1-5\log n/p $, and $1-5\log n/p $,
respectively. For each model we simulated 200 data sets.

We applied the same six methods as those in simulation I to estimate
the sparse $p $-vectors $\bbeta $. For SIS-SCAD and SIS-DS, we chose
$d =[\frac{3}{2}n/\log n]$, $[\frac{3}{2}n/\log n]$, and $[n/\log
n]$, respectively, and for the last two methods, we chose $d=n-1$
and $d' =[\frac{3}{2}n/\log n]$, $[\frac{3}{2}n/\log n]$, and
$[n/\log n]$, respectively. The simulation results are similarly
summarized in Figure 6 (based on 500 simulations) and Table 2.
Similar conclusions as those from simulation I can be drawn. As in
simulation I, we did not include the Gauss-Dantzig selector for
simplicity. It is interesting to observe that in the first setting
here, the Lasso gives large models and its estimation errors are
noticeable compare to the norm of the true coefficient vector
$\bbeta$.

\subsubsection{Leukemia data analysis}
We also applied SIS to select features for the classification of a
Leukemia data set. The Leukemia data from high-density Affymetrix
oligonucleotide arrays were previously analyzed in Golub \textit{et
al.} (1999) and are available at
\verb=http://www.broad.mit.edu/cgi-bin/cancer/datasets.cgi=. There
are 7129 genes and 72 samples from two classes: 47 in class ALL
(acute lymphocytic leukemia) and 25 in class AML (acute mylogenous
leukemia). Among those 72 samples, 38 (27 in class ALL and 11 in
class AML) of them were set as the training sample and the remaining
34 (20 in class ALL and 14 in class AML) of them were set to be the
test sample.

We used two methods SIS-SCAD-LD and SIS-SCAD-NB that will be
introduced below to carry out the classification. For each method,
we first applied SIS to select $d =[2n/\log n]$ genes with $n=38$
the training sample size chosen above and then used the SCAD to get
a family of models indexed by the regularization parameter $\lambda
$. Here, we should point out that our classification results are not
very sensitive to the choice of $d $ as long as it is not too small.
There are certainly many ways to tune the regularization parameter
$\lambda $. For simplicity, we chose a $\lambda $ that produces a
model with size equal to the optimal number of features determined
by the Features Annealed Independence Rules (FAIR) procedure in Fan
and Fan (2008). 16 genes were picked up by their approach. Now we
selected 16 genes and got a linear model with size 16 by using
SIS-SCAD. Finally, the SIS-SCAD-LD method directly used the above
linear discrimination rule to do classification, and the SIS-SCAD-NB
method applied the naive Bayes (NB) rule to the resulted
16-dimensional feature space.

The classification results of the SIS-SCAD-LD, SIS-SCAD-NB, and
nearest shrunken centroids method in Tibshirani \textit{et al.}
(2002) are shown in Table 3. The results of the nearest shrunken
centroids method were extracted from Tibshirani \textit{et al.}
(2002). The SIS-SCAD-LD and SIS-SCAD-NB both chose 16 genes and made
1 test error with training errors 0 and 4, respectively, while the
nearest shrunken centroids method picked up 21 genes and made 1
training error and 2 test errors.

\begin{table}[tph]
\begin{center}
Table 3: Classification errors on the Leukemia data set\vspace{0.2
in}\\\small
\begin{tabular}{@{}c  c  c  c @{}}
\hline
Method & Training error & Test error &  Number of genes  \\
\hline
SIS-SCAD-LD  & 0/38 & 1/34  & 16\\
\hline
SIS-SCAD-NB  & 4/38 & 1/34  & 16\\
\hline
Nearest shrunken centroids   & 1/38 & 2/34  & 21\\
\hline
\end{tabular}
\end{center}
\end{table}

\section{Extensions of SIS}
Like modeling building in linear regression, there are many
variations in the implementation of correlation learning. This
section discusses some extensions of SIS to enhance its
methodological power.  In particular, an iterative SIS (ISIS) is
proposed to overcome some weak points of SIS.  The methodological
power of ISIS is illustrated by three simulated examples.

\subsection{Some extensions of correlation learning}
The key idea of SIS is to apply a single componentwise regression.
Three potential issues, however, might arise with this approach.
First, some unimportant predictors that are highly correlated with
the important predictors can have higher priority to be selected by
SIS than other important predictors that are relatively weakly
related to the response. Second, an important predictor that is
marginally uncorrelated but jointly correlated with the response can
not be picked by SIS and thus will not enter the estimated model.
Third, the issue of collinearity between predictors adds difficulty
to the problem of variable selection. These three issues will be
addressed in the extensions of SIS below, which allow us to use more
fully the joint information of the covariates rather than just the
marginal information in variable selection.

\subsubsection{ISIS: An iterative correlation learning}
It will be shown that when the model assumptions are satisfied,
which excludes basically the three aforementioned problems,  SIS can
accurately reduce the dimensionality from ultra high to a moderate
scale, say, below sample size. But when those assumptions fail, it
could happen that SIS would miss some important predictors. To
overcome this problem, we propose below an ISIS to enhance the
methodological power. It is an iterative applications of the SIS
approach to variable selection.   The essence is to iteratively
apply a large-scale variable screening followed by a moderate-scale
careful variable selection.

The ISIS works as follows. In the first step, we select a subset of
$k_1$ variables $\mathcal{A}_1=\{X_{i_1},\cdots,X_{i_{k_1}}\}$ using
an SIS based model selection method such as the SIS-SCAD or
SIS-Lasso. These variables were selected, using SCAD or Lasso, based
on the joint information of $[n/\log n]$ variables that survive
after the correlation learning. Then we have an $n$-vector of the
residuals from regressing the response $Y$ over
$X_{i_1},\cdots,X_{i_{k_1}}$. In the next step, we treat those
residuals as the new responses and apply the same method as in the
previous step to the remaining $p -k_1$ variables, which results in
a subset of $k_2$ variables
$\mathcal{A}_2=\{X_{j_1},\cdots,X_{j_{k_2}}\}$. We remark that
fitting the residuals from the previous step on $\{X_1,\cdots,X_{p
}\}\setminus\mathcal{A}_1$ can significantly weaken the priority of
those unimportant variables that are highly correlated with the
response through their associations with
$X_{i_1},\cdots,X_{i_{k_1}}$, since the residuals are uncorrelated
with those selected variables in ${\cal A}_1$. This helps solving
the first issue. It also makes those important predictors that are
missed in the previous step possible to survive, which addresses the
second issue above.   In fact, after variables in ${\cal A}_1$
entering into the model, those that are marginally weakly correlated
with $Y$ purely due to the presence of variables in ${\cal A}_1$
should now be correlated with the residuals.  We can keep on doing
this until we get $\ell$ disjoint subsets
$\mathcal{A}_1,\cdots,\mathcal{A}_\ell$ whose union
$\mathcal{A}=\cup_{i=1}^\ell\mathcal{A}_i$ has a size $d$, which is
less than $n$.  In practical implementation, we can choose, for
example, the largest $l$ such that $|\mathcal{A}| < n$.  From the
selected features in $\mathcal{A}$, we can choose the features using
a moderate scale method such as SCAD, Lasso or Dantzig.

For the problem of ultra-high dimensional variable selection, we now
have the ISIS based model selection methods which are extensions of
SIS based model selection methods. Applying a moderate dimensional
method such as the SCAD, Dantzig selector, Lasso, or adaptive Lasso
to $\mathcal{A}$ will produce a model that is very close to the true
sparse model $\mathcal{M}_{*}$.  The idea of ISIS is somewhat
related to the boosting algorithm (Freund and Schapire, 1997).  In
particular, if the SIS is used to select only one variable at each
iteration, i.e., $|\mathcal{A}_i| = 1$, the ISIS is equivalent to a
form of matching pursuit or a greedy algorithm for variable
selection (Barron, {\em et al.}, 2008).

\subsubsection{Grouping and transformation of the input variables}

Grouping the input variables is often used in various problems. For
instance, we can divide the pool of $p $ variables into disjoint
groups each with 5 variables. The idea of variable screening via SIS
can be applied to select a small number of groups. In this way there
is less chance of missing the important variables by taking
advantage of the joint information among the predictors. Therefore a
more reliable model can be constructed.

A notorious difficulty of variable selection lies in the
collinearity between the covariates. Effective ways to rule out
those unimportant variables that are highly correlated with the
important ones are being sought after. A good idea is to transform
the input variables. Two possible ways stand out in this regard. One
is subject related transformation and the other is statistical
transformation.

Subject related transformation is a useful tool. In some cases, a
simple linear transformation of the input variables can help weaken
correlation among the covariates. For example, in somatotype studies
the common sense tells us that predictors such as the weights $w_1,\
w_2$ and $w_3$ at 2, 9 and 18 years are positively correlated. We
could directly use $w_1,\ w_2$ and $w_3$ as the input variables in a
linear regression model, but a better way of model selection in this
case is to use less correlated predictors such as $(w_1, w_2-w_1,
w_3-w_2)^T$, which is a linear transformation of $(w_1, w_2, w_3)^T$
that specifies the changes of the weights instead of the weights
themselves. Another important example is the financial time series
such as the prices of the stocks or interest rates. Differencing can
significantly weaken the correlation among those variables.

Methods of statistical transformation include an application of a
clustering algorithm such as the hierarchical clustering or $k$-mean
algorithm using the correlation metrics to first group variables
into highly correlated groups and then apply the sparse principal
components analysis (PCA) to construct weakly correlated predictors.
Now those weakly correlated predictors from each group can be
regarded as the new covariates and an SIS based model selection
method can be employed to select them.

The statistical techniques we introduced above can help identify the
important features and thus improve the effectiveness of the vanilla
SIS based model selection strategy. Introduction of nonlinear terms
and transformation of variables can also be used to reduced the
modeling biases of linear model. Ravikumar {\em et al.} (2007)
introduced sparse additive models (SpAM) to deal with nonlinear
feature selection.

\subsection{Numerical evidence}
To study the performance of the ISIS proposed above, we now present
three simulated examples. The aim is to examine the extent to which
ISIS can improve SIS in the situation where the conditions of SIS
fail.  We evaluate the methods by counting the frequencies that the
selected models include all the variables in the true model, namely
the ability of correctly screening unimportant variables.

\subsubsection{Simulated example I}
For the first simulated example, we used a linear model
\[ Y=5X_1+5X_2+5X_3+\varepsilon, \]
where $X_1,\cdots,X_p$ are $p$ predictors and $\varepsilon\sim
N(0,1)$ is a noise that is independent of the predictors. In the
simulation, a sample of $(X_1,\cdots,X_p)$ with size $n$ was drawn
from a multivariate normal distribution $N(0,\Sigma)$ whose
covariance matrix $\Sigma=(\sigma_{ij})_{p\times p}$ has entries
$\sigma_{ii}=1$, $i=1,\cdots,p$ and $\sigma_{ij}=\rho$, $i\neq j$.
We considered 20 such models characterized by $(p,n,\rho)$ with
$p=100$, $1000$, $n=20$, $50$, $70$, and $\rho=0$, $0.1$, $0.5$,
$0.9$, respectively, and for each model we simulated 200 data sets.

For each model, we applied SIS and the ISIS to select $n$ variables
and tested their accuracy of including the true model
$\{X_1,X_2,X_3\}$. For the ISIS, the SIS-SCAD with $d=[n/\log n]$
was used at each step and we kept on collecting variables in those
disjoint $\mathcal{A}_j$'s until we got $n$ variables (if there were
more variables than needed in the final step, we only included those
with the largest absolute coefficients). In Table 4, we report the
percentages of SIS, Lasso and ISIS that include the true model.  All
of these three methods select $n-1$-variables, in order to make fair
comparisons. It is clear that the collinearity (large value of
$\rho$) and high-dimensionality deteriorate the performance of SIS
and Lasso, and Lasso outperforms SIS somewhat. However, when the
sample size is 50 or more, the difference in performance is very
small, but SIS has much less computational cost. On the other hand,
ISIS improves dramatically the performance of this simple SIS and
Lasso. Indeed, in this simulation, ISIS always picks all true
variables. It can even have much less computational cost than Lasso
when Lasso is used in the implementation of ISIS.

\begin{table}[tph]
\begin{center}
Table 4: Results of simulated example I: Accuracy of SIS, Lasso and
ISIS \\ in including the true model $\{X_1,X_2,X_3\}$ \vspace{0.2
in}\\\small
\begin{tabular}
[c]{c|cccccc}\hline
$p$ & $n$ &       & \quad$\rho=0$\quad & $\rho=0.1$ & $\rho=0.5$ & $\rho=0.9$ \\
\hline%
    &    & \SIS   & .755    &  .855    &  .690  & .670 \\
    & 20 & \Lasso & .970    &  .990    &  .985  & .870 \\
100 &    & \ISIS  &  1      &   1      &    1   &  1 \\\cline{2-7}
    &    & \SIS   &  1      &   1      &    1   &  1   \\
    & 50 & \Lasso &  1      &   1      &    1   &  1   \\
    &    & \ISIS  &  1      &   1      &    1   &  1 \\
    \hline%
    &    & \SIS   & .205    & .255     &  .145  & .085 \\
    & 20 & \Lasso & .340    & .555     &  .556  & .220 \\
    &    & \ISIS  &  1      &   1      &    1   &  1 \\\cline{2-7}

    &    & \SIS   & .990    & .960     & .870   & .860 \\
1000& 50 & \Lasso &  1      &   1      &    1   &  1   \\
    &    & \ISIS  &  1      &   1      &    1   &  1 \\\cline{2-7}
    &    & \SIS   &  1      & .995     &  .97   &  .97 \\
    & 70 & \Lasso &  1      &   1      &    1   &  1 \\
    &    & \ISIS  &  1      &   1      &    1   &  1 \\  \hline%
\end{tabular}
\end{center}
\end{table}

\subsubsection{Simulated example II}
For the second simulated example, we used the same setup as in
example I except that $\rho$ was fixed to be 0.5 for simplicity. In
addition, we added a fourth variable $X_4$ to the model and the
linear model is now
\[ Y=5X_1+5X_2+5X_3-15\sqrt{\rho}X_4+\varepsilon, \]
where $X_4\sim N(0,1)$ and has correlation $\sqrt{\rho}$ with all
the other $p-1$ variables. The way $X_4$ was introduced is to make
it uncorrelated with the response $Y$.  Therefore, the SIS can not
pick up the true model except by chance.

Again we simulated 200 data sets for each model. In Table 5, we
report the percentages of SIS, Lasso and ISIS that include the true model
of four variables.  In this simulation example, SIS performs somewhat better
than Lasso in variable screening, and ISIS outperforms
significantly the simple SIS and Lasso.   In this simulation it always picks
all true variables.  This demonstrates that ISIS can effectively
handle the second problem mentioned at the beginning of Section 4.1.

\begin{table}[tph]
\begin{center}
Table 5: Results of simulated example II: Accuracy of SIS, Lasso and
ISIS\\ in including the true model $\{X_1,X_2,X_3,X_4\}$ \vspace{0.2
in}\\\small
\begin{tabular}
[c]{cccccc}\hline
$p$  & $\rho=0.5$ &     & \quad$n=20$\quad & $n=50$ & $n=70$ \\
\hline%
     &     & \SIS    &  .025  &  .490   &    .740 \\
100  &     & \Lasso  &  .000  &  .360   &    .915 \\
     &     & \ISIS   &   1    &   1     &     1    \\
\hline%
     &     & \SIS    &  .000  &  .000   &    .000 \\
1000 &     & \Lasso  &  .000  &  .000   &    .000 \\
     &     & \ISIS   &   1    &   1     &     1    \\ \hline
\end{tabular}
\end{center}
\end{table}

\subsubsection{Simulated example III}
For the third simulated example, we used the same setup as in
example II except that we added a fifth variable $X_5$ to the model
and the linear model is now
\[ Y=5X_1+5X_2+5X_3-15\sqrt{\rho}X_4+X_5+\varepsilon, \]
where $X_5\sim N(0,1)$ and is uncorrelated with all the other $p-1$
variables. Again $X_4$ is uncorrelated with the response $Y$. The
way $X_5$ was introduced is to make it have a very small correlation
with the response and in fact the variable $X_5$ has the same
proportion of contribution to the response as the noise
$\varepsilon$ does.  For this particular example, $X_5$ has weaker
marginal correlation with $Y$ than $X_6, \cdots, X_p$ and hence has
a lower priority to be selected by SIS.

For each model we simulated 200 data sets. In Table 6, we report the
accuracy in percentage of SIS, Lasso and ISIS in including the true
model. It is clear to see that the ISIS can improve significantly
over the simple SIS and Lasso and always picks all true
variables.  This shows again that the ISIS is able to pick up two
difficult variables $X_4$ and $X_5$, which addresses simultaneously
the second and third problem at the beginning of Section 4.

\begin{table}[tph]
\begin{center}
Table 6: Results of simulated example III: Accuracy of SIS, Lasso and
ISIS\\ in including the true model $\{X_1,X_2,X_3,X_4,X_5\}$
\vspace{0.2 in}\\\small
\begin{tabular}
[c]{cccccc}\hline
$p$  & $\rho=0.5$ &        & \quad$n=20$\quad & $n=50$ & $n=70$ \\
\hline%
     &      & \SIS   &  .000  & .285   & .645 \\
100  &      & \Lasso &  .000  & .310   & .890 \\
     &      & \ISIS  &  1     &  1     &  1  \\
\hline%
     &      & \SIS   &  .000  & .000   & .000 \\
1000 &      & \Lasso &  .000  & .000   & .000 \\
     &      & \ISIS  &  1     &  1     &  1  \\
\hline%
\end{tabular}
\end{center}
\end{table}

\subsubsection{Simulations I and II in Section 3.3 revisited}
Now let us go back to the two simulation studies presented in
Section 3.3. For each of them, we applied the technique of ISIS with
SCAD and $d=[n/\log n]$ to select $q=[n/\log n]$ variables. After
that, we estimated the $q$-vector $\bbeta$ by using SCAD. This
method is referred to as ISIS-SCAD. We report in Table 7 the median
of the selected model sizes and median of the estimation errors
$\|\hbbeta -\bbeta \|$ in $\ell_2$-norm. We can see clearly that
ISIS improves over the simple SIS.  The improvements are more
drastic for simulation II in which covariates are more correlated
and the variable selections are more challenging.

\begin{table}[tph]
\begin{center}
Table 7: Simulations I and II in Section 3.3 revisited: Medians of
the selected \\model sizes (upper entry) and the estimation errors
(lower entry) \vspace{0.2 in}\\\small
\begin{tabular}
[c]{ccccc}\hline  & \bf Simulation I & & \multicolumn{2}{c}{\bf
Simulation II}
\\\cline{2-2} \cline{4-5}
$p$ & \bf ISIS-SCAD & & & \bf ISIS-SCAD \\
\hline%
1000 & 13 & & $(s=5)$ & 11\\
\cline{2-2} \cline{5-5}
 & 0.329 & & & 0.223\\
 &  & & $(s=8)$ & 13.5\\
  \cline{5-5}
 &  & & & 0.366\\
20000 & 31 & & & 27\\
\cline{2-2} \cline{5-5}
 & 0.246 & & & 0.315\\\hline
\end{tabular}
\end{center}
\end{table}

\section{Asymptotic analysis}
We introduce an asymptotic framework below and present the sure
screening property for both SIS and ITRRS as well as the consistency
of the SIS based model selection methods SIS-DS and SIS-SCAD.

\subsection{Assumptions}
Recall from (\ref{011}) that $Y=\sum_{i=1}^{p }\beta_{i}X_i+\veps$.
Throughout the paper we let $\mathcal{M}_{*}=\left\{1\leq i\leq p
:\beta_{i}\neq 0\right\}$ be the true sparse model with nonsparsity
size $s =|\mathcal{M}_{*}|$ and define
\begin{equation} \label{117}
\bz=\Sig ^{-1/2}\bx \quad\text{and}\quad \bZ=\bX\Sig ^{-1/2},
\end{equation}
where $\bx=\left(X_1,\cdots,X_{p }\right)\t$ and $\Sig
=\cov\left(\bx\right)$. Clearly, the $n$ rows of the transformed
design matrix $\bZ$ are i.i.d. copies of $\bz$ which now has
covariance matrix $I_{p }$. For simplicity, all the predictors
$X_1,\cdots,X_{p }$ are assumed to be standardized to have mean 0
and standard deviation 1. Note that the design matrix $\bX$ can be
factored into $\bZ\Sig ^{1/2}$. Below we will make assumptions on
$\bZ$ and $\Sig $ separately.

We denote by $\lambda_{\text{max}}\left(\cdot\right)$ and
$\lambda_{\text{min}}\left(\cdot\right)$ the largest and smallest
eigenvalues of a matrix, respectively. For $\bZ$, we are concerned
with a concentration property of its extreme singular values as
follows: \\
\textbf{Concentration Property:} The random matrix $\bZ$ is said to
have the concentration property if there exist some $c,c_1>1$ and
$C_1>0$ such that the following deviation inequality
\begin{equation} \label{086}
P\left(\lambda_{\text{max}}(\widetilde{p}
^{-1}\widetilde{\bZ}\widetilde{\bZ}\t)>c_1\text{ and
}\lambda_{\text{min}}(\widetilde{p}
^{-1}\widetilde{\bZ}\widetilde{\bZ}\t)<1/c_1\right)\leq e^{-C_1n}
\end{equation}
holds for any $n\times \widetilde{p} $ submatrix $\widetilde{\bZ}$
of $\bZ$ with $cn<\widetilde{p} \leq p $. We will call it Property C
for short. Property C amounts to a distributional constraint on
$\bz$. Intuitively, it means that with large probability the $n$
nonzero singular values of the $n\times \widetilde{p} $ matrix
$\widetilde{\bZ}$ are of the same order, which is reasonable since
$\widetilde{p} ^{-1}\widetilde{\bZ}\widetilde{\bZ}\t$ will approach
$I_n$ as $\widetilde{p} \rightarrow\infty$: the larger the
$\widetilde{p} $, the closer to $I_n$. It relies on the random
matrix theory (RMT) to derive the deviation inequality in
(\ref{086}). In particular, Property C holds when $\bx$ has a $p
$-variate Gaussian distribution (see Appendix A.7). We conjecture
that it should be shared by a wide class of spherically symmetric
distributions. For studies on the extreme eigenvalues and limiting
spectral distributions, see, e.g., Silverstein (1985), Bai and Yin
(1993), Bai (1999), Johnstone (2001), and Ledoux (2001, 2005).

Some of the assumptions below are purely technical and only serve to
provide theoretical understanding of the newly proposed methodology.
We have no intent to make our assumptions the weakest possible.

\smallskip

\textbf{Condition 1.} $p >n$ and  $ \log p =O(n^\xi)$ for some {$\xi
\in (0, 1-2 \kappa)$, where $\kappa$ is given by Condition 3}.

\smallskip

\textbf{Condition 2.} $\bz$ has a spherically symmetric distribution
and Property C. Also, $\veps\sim\mathcal{N}(0,\sigma^2)$ for some
$\sigma>0$.

\smallskip

\textbf{Condition 3.} $\var\left(Y\right)=O(1)$ and for some
$\kappa\geq0$ and $c_2,c_3>0$, \[
\min_{i\in\mathcal{M}_{*}}\left|\beta_{i}\right|\geq
\frac{c_2}{n^\kappa}\quad\text{and}\quad\min_{i\in\mathcal{M}_{*}}|\cov(\beta_{i}^{-1}Y,X_i)|\geq
c_3. \]

As seen later, $\kappa$ controls the rate of probability error in
recovering the true sparse model. Although $b
=\min_{i\in\mathcal{M}_{*}}|\cov(\beta_{i}^{-1}Y,X_i)|$ is assumed
here to be bounded away from zero, our asymptotic study applies as
well to the case where $b $ tends to zero as $n\rightarrow\infty$.
In particular, when the variables in ${\cal M}_*$ are uncorrelated,
$b = 1$. This condition rules out the situation in which an
important variable is marginally uncorrelated with $Y$, but jointly
correlated with $Y$.

\smallskip

\textbf{Condition 4.} There exist some $\tau\geq0$ and $c_4>0$ such
that
\[ \lambda_{\text{max}}\left(\Sig \right)\leq c_4n^\tau.
\]
This condition rules out the case of strong collinearity.

The largest eigenvalue of the population covariance matrix $\Sig $
is allowed to diverge as $n$ grows. When there are many predictors,
it is often the case that their covariance matrix is block diagonal
or nearly block diagonal under a suitable permutation of the
variables. Therefore $\lambda_{\text{max}}\left(\Sig \right)$
usually does not grow too fast with $n$. In addition, Condition 4
holds for the covariance matrix of a stationary time series (see
Bickel and Levina, 2004, 2008). See also Grenander and Szeg\"{o}
(1984) for more details on the characterization of extreme
eigenvalues of the covariance matrix of a stationary process in
terms of its spectral density.

\subsection{Sure screening property}
Analyzing the $p$-vector $\bomega$ in (\ref{120}) when $p >n$ is
essentially difficult. The approach we took is to first study the
specific case with $\Sig =I_{p }$ and then relate the general case
to the specific case.

\begin{theorem}
{\em (Accuracy of SIS).} Under Conditions 1--4, if $2\kappa+\tau<1$
then there exists some $\theta<1-2\kappa-\tau$ such that when
$\gamma \sim cn^{-\theta}$ with $c>0$,  we have for some $C>0$,
\[ P\left(\mathcal{M}_{*}\subset\mathcal{M}_{\gamma }\right)
=1-O(\exp(-Cn^{1-2\kappa}/\log n)).
\]
\end{theorem}

We should point out here that $s \leq[\gamma n]$ is implied by our
assumptions as demonstrated in the technical proof. The above
theorem shows that SIS has the sure screening property and can
reduce from exponentially growing dimension $p $ down to a
relatively large scale $d =[\gamma n]= O(n^{1-\theta})<n$ for some
$\theta>0$, where the reduced model $\mathcal{M}=\mathcal{M}_{\gamma
}$ still contains all the variables in the true model with an
overwhelming probability. In particular, we can choose the submodel
size $d$ to be $n-1$ or $n/\log n$ for SIS if Conditions 1-4 are
satisfied.

Another interpretation of Theorem 1 is that it requires the model
size $d = [\gamma n] = n^{\theta^*}$ with $\theta^* > 2 \kappa +
\tau$ in order to have the sure screening property.  The weaker the
signal, the larger the $\kappa$ and hence the larger the required
model size. Similarly, the more severe the collinearity, the larger
the $\tau$ and the larger the required model size. In this sense,
the restriction that $2\kappa + \tau < 1$ is not needed, but $\kappa
< 1/2$ is needed since we can not detect signals that of smaller
order than root-$n$ consistent.  In the former case, there is no
guarantee that $\theta^*$ can be taken to be smaller than one.

The proof of Theorem 1 depends on the iterative application of the
following theorem, which demonstrates the accuracy of each step
of ITRRS.  We first describe the result of the first step of ITRRS.
It shows that as long as the ridge parameter $\lambda$ is large
enough and the percentage of remaining variables $\delta$ is large
enough, the sure screening property is ensured with overwhelming
probability.

\begin{theorem}
{\em (Asymptotic sure screening).} Under Conditions 1--4, if
$2\kappa+\tau<1$, $\lambda (p ^{3/2}n)^{-1}\rightarrow \infty$, and
$\delta n^{1-2\kappa-\tau}\rightarrow\infty$ as
$n\rightarrow\infty$, then we have for some $C>0$,
\[
 P\left(\mathcal{M}_{*}\subset\mathcal{M}^{1}_{\delta, \lambda}\right)
  =1-O(\exp(-Cn^{1-2\kappa}/\log
n)).
\]
\end{theorem}

The above theorem reveals that when the tuning parameters are chosen
appropriately, with an overwhelming probability the submodel
$\mathcal{M}^{1}_{\delta, \lambda}$ will contain the true model
$\mathcal{M}_{*}$ and its size is an order $n^\theta$ (for some
$\theta > 0$) lower than the original one. This property stimulated
us to propose ITRRS.

\begin{theorem}
{\em (Accuracy of ITRRS).} Let the assumptions of Theorem 2 be
satisfied. If $\delta n^\theta\rightarrow \infty$ as
$n\rightarrow\infty$ for some $\theta<1-2\kappa-\tau$, then
successive applications of the procedure in (\ref{039}) for $k$
times results in a submodel $\mathcal{M}_{\delta, \lambda}$ with
size $d =[\delta ^{k }p ] < n$ such that for some $C>0$,
\[
P\left(\mathcal{M}_{*}\subset\mathcal{M}_{\delta, \lambda}
\right)=1-O(\exp(-Cn^{1-2\kappa}/\log n)).
\]
\end{theorem}

Theorem 3 follows from iterative application of Theorem 2 $k$ times,
where $k$ is the first integer such that $[\delta^k p] < n$.  This
implies that $k = O(\log p / \log n) = O(n^\xi)$.  Therefore, the
accumulated error probability, from the union bound, is still of
exponentially small with a possibility of a different constant $C$.

ITRRS has now been shown to possess the sure screening property. As
mentioned before, SIS is a specific case of ITRRS with an infinite
regularization parameter and hence enjoys also the sure screening
property.

Note that the number of steps in ITRRS depends on the choice of
$\delta\in(0,1)$. In particular, $\delta$ can not be too small, or
equivalently, the number of iteration steps in ITRRS can not be too
large, due to the accumulation of the probability errors of missing
some important variables over the iterations. In particular, the
stepwise deletion method which deletes one variable each time in
ITRRS might not work since it requires $p-d$ steps of iterations,
which may exceed the error bound in Theorem 2.

\subsection{Consistency of SIS-DS and SIS-SCAD}
To study the property of the Dantzig selector, Candes and Tao (2007)
introduce the notion of uniform uncertainty principle (UUP) on
deterministic design matrices which essentially states that the
design matrix obeys a ``restricted isometry hypothesis."
Specifically, let $\bA$ be an $n\times d $ deterministic design
matrix and for any subset $T\subset\{1,\cdots,d\}$. Denote by
$\bA_T$ the $n\times |T|$ submatrix of $\bA$ obtained by extracting
its columns corresponding to the indices in $T$. For any positive
integer $S\leq d $, the $S$-restricted isometry constant
$\delta_S=\delta_S(\bA)$ of $\bA$ is defined to be the smallest
quantity such that
\[ \left(1-\delta_S\right)\left\|\bv\right\|^2\leq\left\|\bA_T\bv\right\|^2
\leq\left(1+\delta_S\right)\left\|\bv\right\|^2 \] holds for all
subsets $T$ with $|T|\leq S$ and $\bv\in\mathbf{R}^{|T|}$. For any
pair of positive integers $S,S'$ with $S+S'\leq d $, the
$S,S'$-restricted orthogonality constant
$\theta_{S,S'}=\theta_{S,S'}(\bA)$ of $\bA$ is defined to be the
smallest quantity such that
\[ \left|\left\langle\bA_T\bv,\bA_{T'}\bv'\right\rangle\right|\leq\theta_{S,S'}
\left\|\bv\right\|\left\|\bv'\right\| \] holds for all disjoint
subsets $T,T'$ of cardinalities $|T|\leq S$ and $|T'|\leq S'$,
$\bv\in\mathbf{R}^{|T|}$, and $\bv'\in\mathbf{R}^{|T'|}$.

The following theorem is obtained by the sure screening property of
SIS in Theorem 1 along with Theorem 1.1 in Candes and Tao (2007),
where $\bveps\sim\mathcal{N}(\bzero,\sigma^2I )$ for some
$\sigma>0$. To avoid the selection bias in the prescreening step, we
can split the sample into two halves:  the first half is used to
screen variables and the second half is used to construct the
Dantzig estimator. The same technique applies to SCAD, but we avoid
this step of detail for simplicity of presentation.

\begin{theorem}
{\em (Consistency of SIS-DS).} Assume with large probability,
$\delta_{2s }(\bX_\mathcal{M})+\theta_{s ,2s }(\bX_\mathcal{M})\leq
t<1$ and choose $\lambda_{d }=\sqrt{2\log d }$ in (\ref{002}). Then
with large probability, we have
\[ \left\|\hbbeta_{\text{\em DS}}-\bbeta \right\|^2\leq C\left(\log d \right)s \sigma^2, \]
where $C=32/\left(1-t\right)^2$ and $s $ is the number of nozero
components of $\bbeta$.
\end{theorem}

This theorem shows that SIS-DS, i.e., SIS followed by the Dantzig
selector, can now achieve the ideal risk up to a factor of $\log d $
with $d<n$, rather than the original $\log p$.

Now let us look at SIS-SCAD, that is, SIS followed by the SCAD. For
simplicity, a common regularization parameter $\lambda $ is used for
the SCAD penalty function. Let
$\hbbeta_{\text{SCAD}}=\left(\hbeta_1,\cdots,\hbeta_{d }\right)\t$
be a minimizer of the SCAD-PLS in (\ref{111}). The following theorem
is obtained by the sure screening property of SIS in Theorem 1 along
with Theorems 1 and 2 in Fan and Peng (2004).

\begin{theorem}
{\em (Oracle properties of SIS-SCAD).} If $d =o(n^{1/3})$ and the
assumptions of Theorem 2 in Fan and Peng (2004) be satisfied, then,
with probability tending to one, the SCAD-PLS estimator
$\hbbeta_{\text{\em SCAD}}$ satisfies: (i) $\hbeta_i=0$ for any $i
\not \in\mathcal{M}_{*}$; (ii) the components of $\hbbeta_{\text{\em
SCAD}}$ in $\mathcal{M}_{*}$ perform as well as if the true model
$\mathcal{M}_{*}$ were known.
\end{theorem}

The SIS-SCAD has been shown to enjoy the oracle properties.

\section{Concluding remarks}
This paper studies the problem of high dimensional variable
selection for the linear model. The concept of sure screening is
introduced and a sure screening method based on correlation learning
that we call the Sure Independence Screening (SIS) is proposed. The
SIS has been shown to be capable of reducing from exponentially
growing dimensionality to below sample size accurately. It speeds up
variable selection dramatically and can also improve the estimation
accuracy when dimensionality is ultra high. SIS combined with
well-developed variable selection techniques including the SCAD,
Dantzig selector, Lasso, and adaptive Lasso provides a powerful tool
for high dimensional variable selection.  The tuning parameter $d$
can be taken as  $d=[n/\log n]$ or $d = n-1$, depending on which
model selector is used in the second stage.  For non-concave
penalized least-squares (\ref{jf1}), when one directly applies the
LLA algorithm to the original problem with $d = p$, one needs
initial values that are not readily available.  SIS provides a
method that makes this feasible by screening many variables and
furnishing the corresponding coefficients with zero.  The initial
value in (\ref{jf1}) can be taken as the OLS estimate if $d=[n/\log
n]$ and zero [corresponding to $w_j^{(0)} \equiv p_\lambda'(0+)$]
when $d = n-1$, which is LASSO.

Some extensions of SIS have also been discussed. In particular, an
iterative SIS (ISIS) is proposed to enhance the finite sample
performance of SIS, particularly in the situations where the
technical conditions fail. This raises a challenging question: to
what extent does ISIS relax the conditions for SIS to have the sure
screening property? An iteratively thresholded ridge regression
screener (ITRRS) has been introduced to better understand the
rationale of SIS and serves as a technical device for proving the
sure screening property. As a by-product, it is demonstrated that
the stepwise deletion method may have no sure screening property
when the dimensionality is of an exponential order. This raises
another interesting question if the sure screening property holds
for a greedy algorithm such as the stepwise addition or matching
pursuit and how large the selected model has to be if it does.

The paper leaves open the problem of extending the SIS and ISIS
introduced for the linear models to the family of generalized linear
models (GLM) and other general loss functions such as the hinge loss
and the loss associated with the support vector machine (SVM).
Questions including how to define associated residuals to extend
ISIS and whether the sure screening property continues to hold
naturally arise. The paper focuses only on random designs which
commonly appear in statistical problems, whereas for many problems
in fields such as image analysis and signal processing the design
matrices are often deterministic. It remains open how to impose a
set of conditions that ensure the sure screening property.  It also
remains open if the sure screening property can be extended to the
sparse additive model in nonparametric learning as studied by
Ravikumar {\em et al.} (2007). These questions are beyond the scope
of the current paper and are interesting topics for future research.

\appendix

\section{Appendix}
Hereafter we use both $C$ and $c$ to denote generic positive constants for notational convenience.

\subsection{Proof of Theorem 1}
Motivated by the results in Theorems 2 and 3, the idea is to
successively apply dimensionality reduction in a way described in
(\ref{122}) below. To enhance the readability, we split the whole
proof into two mains steps and multiple substeps.

{\bf Step 1.}\quad Let $\delta \in\left(0,1\right)$. Similarly to
(\ref{039}), we define a submodel
\begin{equation} \label{122}
\widetilde{\mathcal{M}}^1_{ \delta }=\left\{1\leq i\leq p
:|\omega_i|\text{ is among the first $\left[\delta p \right]$
largest of all}\right\}.
\end{equation}
We aim to show that if $\delta  \to 0$ in such a way that $\delta
n^{1-2\kappa-\tau}\rightarrow\infty$ as $n\rightarrow\infty$, we
have for some $C>0$,
\begin{equation} \label{123}
P\left(\mathcal{M}_{ *}\subset\widetilde{\mathcal{M}}^1_{ \delta
}\right)= 1-O(\exp(-Cn^{1-2\kappa}/\log n)).
\end{equation}

The main idea is to relate the general case to the specific case
with $\Sig =I_{p }$, which is separately studied in Sections
A.4--A.6 below. A key ingredient is the representation (\ref{089})
below of the $p \times p $ random matrix $\bX\t\bX$. Throughout, let
$\bS = \left(\bZ\t\bZ\right)^+\bZ\t\bZ$ and
$\be_i=\left(0,\cdots,1,\cdots,0\right)\t$ be a unit vector in
$\mathbf{R}^{p }$ with the $i$-th entry 1 and 0 elsewhere,
$i=1,\cdots,p $.

Since $\bX=\bZ\Sig ^{1/2}$, it follows from (\ref{088}) that
\begin{equation} \label{089}
\bX\t\bX=p \Sig ^{1/2}\widetilde{\bU}\t\diag\left(\mu_1,\cdots,\mu_n
\right)\widetilde{\bU}\Sig ^{1/2},
\end{equation}
where $\mu_1,\cdots,\mu_n$ are $n$ eigenvalues of $p ^{-1}\bZ\bZ\t$,
$\widetilde{\bU}=\left(I_n ,\bzero\right)_{n\times p }\bU$, and
$\bU$ is uniformly distributed on the orthogonal group
$\mathcal{O}(p )$. By  (\ref{011}) and (\ref{120}), we have
\begin{equation} \label{041}
\bomega=\bX\t\bX\bbeta +\bX\t\bveps\heq\bxi+\bleta.
\end{equation}
We will study the above two random vectors $\bxi$ and $\bleta$
separately.

{\em Step 1.1}.\quad First, we consider term
$\bxi=\left(\xi_1,\cdots,\xi_{p }\right)\t=\bX\t\bX\bbeta $.

{\em Step 1.1.1. Bounding $\left\|\bxi\right\|$ from above}.\quad It
is obvious that
\[ \diag\left(
\mu_1^2,\cdots,\mu_n ^2\right)\leq\left[\lambda_{\text{max}}(p
^{-1}\bZ\bZ\t)\right]^2I_n
\] and $\widetilde{\bU}\Sig \widetilde{\bU}\t \leq
\lambda_{\text{max}}(\Sig )I_n $. These and (\ref{089}) lead to
\begin{align} \label{093}
\left\|\bxi\right\|^2 \leq p ^2\lambda_{\text{max}}(\Sig
)\left[\lambda_{\text{max}}(p ^{-1}\bZ\bZ\t)\right]^2 \bbeta \t\Sig
^{1/2}\widetilde{\bU}\t\widetilde{\bU}\Sig ^{1/2}\bbeta .
\end{align}

Let $Q\in\mathcal{O}(p )$ such that $\Sig ^{1/2}\bbeta =\left\|\Sig
^{1/2}\bbeta \right\|Q\be_1$. Then, it follows from Lemma 1 that
\begin{align*}
\bbeta \t\Sig ^{1/2}\widetilde{\bU}\t\widetilde{\bU}\Sig ^{1/2}
\bbeta  =\left\|\Sig ^{1/2}\bbeta \right\|^2 \left\langle Q\t \bS
Q\be_1,\be_1\right\rangle \deq\left\|\Sig ^{1/2}\bbeta \right\|^2
\left\langle \bS \be_1,\be_1\right\rangle,
\end{align*}
where we use the symbol $\deq$ to denote being identical in
distribution for brevity. By Condition 3, $\left\|\Sig ^{1/2}\bbeta
\right\|^2=\bbeta \t\Sig \bbeta \leq\var\left(Y\right)=O(1)$, and
thus by Lemma 4, we have for some $C>0$,
\begin{equation} \label{094}
P\left(\bbeta \t\Sig ^{1/2}\widetilde{\bU}\t\widetilde{\bU}\Sig
^{1/2}\bbeta >O(\frac{n}{p })\right)\leq O(e^{-Cn}).
\end{equation}
Since $\lambda_{\text{max}}\left(\Sig \right)=O(n^\tau)$ and
$P\left(\lambda_{\text{max}}(p ^{-1}\bZ\bZ\t)>c_1\right)\leq
e^{-C_1n}$ by Conditions 2 and 4, (\ref{093}) and (\ref{094}) along
with Bonferroni's inequality yield
\begin{equation} \label{095}
P\left(\left\|\bxi\right\|^2>O(n^{1+\tau}p )\right)\leq O(e^{-Cn}).
\end{equation}

{\em Step 1.1.2. Bounding $\left|\xi_i\right|$, $i\in\mathcal{M}_{
*}$, from below}.\quad This needs a delicate analysis. Now fix an
arbitrary $i\in\mathcal{M}_{ *}$. By (\ref{089}), we have
\[ \xi_i=p \be_i\t\Sig ^{1/2}\widetilde{\bU}\t\diag\left(
\mu_1,\cdots,\mu_n\right)\widetilde{\bU}\Sig ^{1/2}\bbeta . \] Note
that $\left\|\Sig
^{1/2}\be_i\right\|=\sqrt{\var\left(X_i\right)}=1$, $\left\|\Sig
^{1/2}\bbeta \right\|=O(1)$. By Condition 3, there exists some $c>0$
such that
\begin{equation} \label{098}
\left|\left\langle\Sig ^{1/2}\bbeta ,\Sig
^{1/2}\be_i\right\rangle\right| =\left|\beta_{ i}\right|\left|\cov
\left(\beta_{ i}^{-1}Y,X_i\right)\right|\geq c/n^\kappa.
\end{equation}
Thus, there exists $Q\in\mathcal{O}(p )$ such that $\Sig
^{1/2}\be_i=Q\be_1$ and
\[ \Sig ^{1/2}\bbeta =\left\langle\Sig ^{1/2}\bbeta ,\Sig ^{1/2}\be_i
\right\rangle Q\be_1+O(1)Q\be_2. \]
Since $\left(\mu_1,\cdots,\mu_n\right)\t$ is independent of
$\widetilde{\bU}$ by Lemma 1 and the uniform distribution on the
orthogonal group $\mathcal{O}(p )$ is invariant under itself, it
follows that
\begin{align} \label{099}
\xi_i \deq & p \left\langle\Sig ^{1/2}\bbeta ,\Sig
^{1/2}\be_i\right\rangle R_1 + O(p )R_2 \heq\xi_{i,1}+\xi_{i,2},
\end{align}
where $\bR = \left(R_1,R_2,\cdots,R_{p }\right)\t =
\widetilde{\bU}\t\diag\left( \mu_1,\cdots,\mu_n
\right)\widetilde{\bU}\be_1$.  We will examine the above two terms
$\xi_{i,1}$ and $\xi_{i,2}$ separately. Clearly,
\[ R_1\geq\be_1\t
\widetilde{\bU}\t\lambda_{\text{min}}(p ^{-1}\bZ\bZ\t)I_n
\widetilde{\bU}\be_1= \lambda_{\text{min}}(p
^{-1}\bZ\bZ\t)\left\langle \bS \be_1,\be_1\right\rangle,
\]
and thus by Condition 2, Lemma 4, and Bonferroni's inequality, we
have for some $c>0$ and $C>0$,
\[
P\left(R_1<cn/{p }\right)\leq O(e^{-Cn}).
\]
This along with (\ref{098}) gives for some $c>0$,
\begin{equation} \label{104}
P\left(\left|\xi_{i,1}\right|<cn^{1-\kappa}\right)\leq O(e^{-Cn}).
\end{equation}

Similarly to Step 1.1.1, it can be shown that
\begin{equation} \label{101}
P\left( \| \bR \| ^2>O(n/{p })\right)\leq O(e^{-Cn}).
\end{equation}
Since $\left(\mu_1,\cdots,\mu_n\right)\t$ is independent of
$\widetilde{\bU}$ by Lemma 1, the argument in the proof of Lemma 5
applies to show that the distribution of $\tilde{\bR} =
\left(R_2,\cdots,R_{p }\right)\t$ is invariant under the orthogonal
group $\mathcal{O}(p -1)$. Then, it follows that $\tilde{\bR} \deq
\| \tilde{\bR} \| \; \bW /\| \bW\|$, where $\bW = (W_1,\cdots,W_{p
-1})^T \sim \mathcal{N}(0,I_{p-1})$, independent of
$\|\tilde{\bR}\|$. Thus, we have
\begin{equation} \label{102}
R_2\deq \|\tilde{\bR}\| W_1/ \| \bW \| .
\end{equation}
In view of (\ref{101}), (\ref{102}), and $\xi_{i,2}=O(p R_2)$,
applying the argument in the proof of Lemma 5 gives for some $c>0$,
\begin{equation} \label{103}
P\left(\left|\xi_{i,2}\right|>c\sqrt{n}|W| \right) \leq O(e^{-Cn}),
\end{equation}
where $W$ is a $\mathcal{N}(0,1)$-distributed random variable.

Let $ x_n =c \sqrt{2C}n^{1-\kappa}/\sqrt{\log n}$. Then, by the
classical Gaussian tail bound, we have
\begin{align*}
P\left(c \sqrt{n} |W| >x_n \right) \leq\sqrt{2/\pi}\
\frac{\exp\left(-Cn^{1-2\kappa}/\log n\right)}{\sqrt{2C}\
n^{1/2-\kappa}/\sqrt{\log n}} =O(\exp(-Cn^{1-2\kappa}/\log n)),
\end{align*}
which along with (\ref{103})and Bonferroni's inequality shows that
\begin{equation} \label{105}
P\left( \left|\xi_{i,2}\right|>{x_n} \right)\leq  =
O(\exp(-Cn^{1-2\kappa}/\log n)).
\end{equation}
Therefore, by Bonferroni's inequality, combining (\ref{099}),
(\ref{104}), and (\ref{105}) together gives for some $c>0$,
\begin{equation} \label{155}
P\left(\left|\xi_i\right|<cn^{1-\kappa}\right)\leq
O(\exp(-Cn^{1-2\kappa}/\log n)),\quad i\in\mathcal{M}_{ *}.
\end{equation}

{\em Step 1.2}.\quad Then, we examine term
$\bleta=\left(\eta_1,\cdots,\eta_{p }\right)\t=\bX\t\bveps$.

{\em Step 1.2.1. Bounding $\left\|\bleta\right\|$ from above}.\quad
Clearly, we have
\begin{align*}
\bX\bX\t =\bZ\Sig \bZ\t\leq\bZ\lambda_{\text{max}}(\Sig )I_{p }\bZ\t
=p \lambda_{\text{max}}(\Sig )\lambda_{\text{max}}(p
^{-1}\bZ\bZ\t)I_n .
\end{align*}
Then, it follows that
\begin{align} \label{106}
\left\|\bleta\right\|^2&=\bveps\t\bX\bX\t\bveps\leq p
\lambda_{\text{max}}(\Sig )\lambda_{\text{max}}(p
^{-1}\bZ\bZ\t)\left\|\bveps\right\|^2.
\end{align}
From Condition 2, we know that $\veps_1^2/\sigma^2,\cdots,\veps_n
^2/\sigma^2$ are i.i.d. $\chi_1^2$-distributed random variables.
Thus, by (\ref{004}) in Lemma 3, there exist some $c>0$ and $C>0$
such that
\[
P\left(\left\|\bveps\right\|^2>cn\sigma^2\right)\leq e^{-Cn},
\]
which along with (\ref{106}), Conditions 2 and 4, and Bonferroni's
inequality yields
\begin{equation} \label{107}
P\left(\left\|\bleta\right\|^2>O(n^{1+\tau}p )\right)\leq
O(e^{-Cn}).
\end{equation}

{\em Step 1.2.2. Bounding $\left|\eta_i\right|$ from above}.\quad
Given that $\bX=X$, $\bleta=X\t\bveps \sim {\cal N}(\bzero,
\sigma^2X\t X)$. Hence,
$\eta_i|_{\bX=X}\sim\mathcal{N}(0,\var\left(\eta_i|\bX=X\right))$
with
\begin{equation} \label{050}
\var\left(\eta_i|\bX=X\right)=\sigma^2\be_i\t X\t X\be_i.
\end{equation}
Let $\mathcal{E}$ be the event
$\left\{\var\left(\eta_i|\bX\right)\leq cn\right\}$ for some $c>0$.
Then, using the same argument as that in Step 1.1.1, we can easily
show that for some $C>0$,
\begin{equation} \label{108}
P\left(\mathcal{E}^c\right)\leq O(e^{-Cn}).
\end{equation}
On the event $\mathcal{E}$, we have
\begin{equation} \label{058}
P\left(\left|\eta_i\right|>x|\bX\right)\leq P\left(\sqrt{cn}|W|>x
\right) \text{ for any } x>0,
\end{equation}
where $W$ is a $\mathcal{N}(0,1)$-distributed random variable. Thus,
it follows from (\ref{108}) and (\ref{058}) that
\begin{align} \label{059}
P\left(\left|\eta_i\right|>x\right) \leq O(e^{-Cn})+
P\left(\sqrt{cn} {\left|W\right|}>x\right).
\end{align}

Let $ x_n'=\sqrt{2cC}n^{1-\kappa}/\sqrt{\log n}$.   Then, invoking
the classical Gaussian tail bound again,  we have
\begin{align*}
P\left(\sqrt{cn} {\left|W\right|}>x_n'\right)
&=O(\exp(-Cn^{1-2\kappa}/\log n)),
\end{align*}
which along with (\ref{059}) {and Condition 1}  shows that
\begin{equation} \label{109}
P\left({\max_{i}} \left|\eta_i\right|>o(n^{1-\kappa})\right)\leq
{O(p\exp(-Cn^{1-2\kappa}/\log n)) =}  O(\exp(-Cn^{1-2\kappa}/\log
n)).
\end{equation}

{\em Step 1.3}.\quad Finally, we combine the results obtained in
Steps 1.1 and 1.2 together. By Bonferroni's inequality, it follows
from (\ref{041}), (\ref{095}), (\ref{155}), (\ref{107}), and
(\ref{109}) that for some constants $c_1,c_2,C>0$,
\begin{equation} \label{119}
P\left(\min_{i\in\mathcal{M}_{
*}}\left|\omega_i\right|<c_1n^{1-\kappa} \text{ or
}\left\|\bomega\right\|^2>c_2n^{1+\tau}p \right)\leq O(s
\exp(-Cn^{1-2\kappa}/\log n)).
\end{equation}
This shows that with overwhelming probability $1-O(s
\exp(-Cn^{1-2\kappa}/\log n))$, the magnitudes of $\omega_i$,
$i\in\mathcal{M}_{ *}$, are uniformly at least of order
$n^{1-\kappa}$ and more importantly, for some $c>0$,
\begin{equation} \label{124}
\#\left\{1\leq k\leq p
:\left|\omega_k\right|\geq\min_{i\in\mathcal{M}_{
*}}\left|\omega_i\right|\right\} \leq c\frac{n^{1+\tau}p
}{\left(n^{1-\kappa}\right)^2} =\frac{cp }{n^{1-2\kappa-\tau}},
\end{equation}
where $\#\{\cdot\}$ denotes the number of elements in a set.

Now, we are ready to see from (\ref{124}) that if $\delta $
satisfies $\delta n^{1-2\kappa-\tau}\rightarrow\infty$ as
$n\rightarrow\infty$, then (\ref{123}) holds for some constant $C>0$
larger than that in (\ref{119}).

{\bf Step 2.}\quad Fix an arbitrary $r\in(0,1)$ and choose a
shrinking factor $\delta $ of the form $(\frac{n}{p })^{\frac{1}{k
-r}}$, for some integer $k \geq 1$. We successively perform
dimensionality reduction until the number of remaining variables drops to
below sample size $n$:
\begin{itemize}
\item First, carry out the
procedure in (\ref{122}) to the full model
$\widetilde{\mathcal{M}}^0_{ \delta }\heq\left\{1,\cdots,p \right\}$
and get a submodel $\widetilde{\mathcal{M}}^1_{ \delta }$ with size
$[\delta p ]$;

\item Then, apply a similar procedure to the model
$\widetilde{\mathcal{M}}^1_{ \delta }$ and again obtain a submodel
$\widetilde{\mathcal{M}}^2_{ \delta
}\subset\widetilde{\mathcal{M}}^1_{ \delta }$ with size $[\delta ^2p
]$, and so on;

\item Finally, get a submodel
$\widetilde{\mathcal{M}}_{ \delta }\heq\widetilde{\mathcal{M}}^{k
}_{ \delta }$ with size $d =[\delta ^{k }p ]=[\delta ^rn]<n$, where
$[\delta ^{k -1}p ]=[\delta ^{r-1}n]>n$.
\end{itemize}
It is obvious that $\widetilde{\mathcal{M}}_{ \delta }=\mathcal{M}_{
\gamma }$, where $\gamma =\delta ^r<1$.

Now fix an arbitrary $\theta_1\in(0,1-2\kappa-\tau)$ and pick some
$r<1$ very close to 1 such that
$\theta_0=\theta_1/r<1-2\kappa-\tau$. We choose a sequence of
integers $k \geq 1$ in a way such that
\begin{equation} \label{125}
\delta n^{1-2\kappa-\tau}\rightarrow\infty \quad\text{and}\quad
\delta n^{\theta_0}\rightarrow0\quad \text{as }n\rightarrow\infty,
\end{equation}
where $\delta =(\frac{n}{p })^{\frac{1}{k -r}}$. Then, applying the
above scheme of dimensionality reduction results in a submodel
$\widetilde{\mathcal{M}}_{ \delta }=\mathcal{M}_{ \gamma }$, where
$\gamma =\delta ^r$ satisfies
\begin{equation} \label{126}
\gamma n^{r(1-2\kappa-\tau)}\rightarrow\infty \quad\text{and}\quad
\gamma n^{\theta_1}\rightarrow0\quad \text{as } n\rightarrow\infty.
\end{equation}

Before going further, let us make two important observations. First,
for any principal submatrix $\Sig^0 $ of $\Sig $ corresponding to a
subset of variables, Condition 4 ensures that
\[
\lambda_{\text{max}}\left(\Sig^0
\right)\leq\lambda_{\text{max}}\left(\Sig \right)\leq c_4n^\tau.
\]
Second, by definition, Property C in (\ref{086}) holds for any
$n\times \widetilde{p} $ submatrix $\widetilde{\bZ}$ of $\bZ$ with
$cn<\widetilde{p} \leq p $, where $c>1$ is some constant. Thus, the
probability bound in (\ref{123}) is uniform over dimension
$\widetilde{p} \in(cn,p ]$. Therefore, for some $C>0$, by
(\ref{125}) and (\ref{123}) we have in each step $1\leq i\leq k $ of
the above dimensionality reduction,
\[
P\left(\mathcal{M}_{ *}\subset\widetilde{\mathcal{M}}^i_{ \delta }|
\mathcal{M}_{ *}\subset\widetilde{\mathcal{M}}^{i-1}_{ \delta
}\right)= 1-O(\exp(-Cn^{1-2\kappa}/\log n)),
\]
which along with Bonferroni's inequality gives
\begin{equation} \label{127}
P\left(\mathcal{M}_{ *}\subset\mathcal{M}_{ \gamma }\right)= 1-O(k
\exp(-Cn^{1-2\kappa}/\log n)).
\end{equation}

It follows from (\ref{125}) that $k =O(\log p /\log n)$, which is of
order $O(n^\xi/\log n)$ by Condition 1. Thus, a suitable increase of
the constant $C>0$ in (\ref{127}) yields
\[
P\left(\mathcal{M}_{ *}\subset\mathcal{M}_{ \gamma
}\right)=1-O(\exp(-Cn^{1-2\kappa}/\log n)).
\]
Finally, in view of (\ref{126}), the above probability bound holds
for any $\gamma \sim cn^{-\theta}$, with $\theta<1-2\kappa-\tau$ and
$c>0$. This completes the proof.

\subsection{Proof of Theorem 2}
One observes that (\ref{039}) uses only the order of componentwise
magnitudes of $\bomega^{\lambda }$, so it is invariant under
scaling. Therefore, in view of (\ref{096}) we see from Step 1 of the
proof of Theorem 1 that Theorem 2 holds for sufficiently large
regularization parameter $\lambda $.

It remains to specify a lower bound on $\lambda $. Now we rewrite
the $p$-vector $\lambda\bomega^{\lambda }$ as
\[ \lambda\bomega^{\lambda }=\bomega-\left[I_{p }-\left(I_{p }+\lambda ^{-1}
\bX\t\bX\right)^{-1}\right]\bomega. \] Let
$\bzeta=\left(\zeta_1,\cdots,\zeta_{p }\right)\t=\left[I_{p
}-\left(I_{p }+\lambda ^{-1} \bX\t\bX\right)^{-1}\right]\bomega$. It
follows easily from $\bX\t\bX=\Sig ^{1/2}\bZ\t\bZ\Sig ^{1/2}$ that
\[ \lambda_{\text{max}}(\bX\t\bX)\leq
p \lambda_{\text{max}}(p ^{-1}\bZ\bZ\t)\lambda_{\text{max}}(\Sig ),
\]
and thus
\begin{align*}
\left\|\bzeta\right\|^2&\leq\left[\lambda_{\text{max}}\left(I_{p
}-(I_{p }+\lambda ^{-1}
\bX\t\bX)^{-1}\right)\right]^2\left\|\bomega\right\|^2\\
&\leq\left[\lambda_{\text{max}}(\lambda ^{-1}
\bX\t\bX)\right]^2\left\|\bomega\right\|^2\\
&\leq \lambda ^{-2}p ^2\left[\lambda_{\text{max}}(p
^{-1}\bZ\bZ\t)\right]^2 \left[\lambda_{\text{max}}(\Sig
)\right]^2\left\|\bomega\right\|^2,
\end{align*}
which along with (\ref{119}), Conditions 2 and 4, and Bonferroni's
inequality shows that
\[ P\left(\left\|\bzeta\right\|>O(\lambda ^{-1}n^{\frac{1+3\tau}{2}}p ^{3/2})\right)\leq
O(s \exp(-Cn^{1-2\kappa}/\log n)). \] Again, by Bonferroni's
inequality and (\ref{119}), any $\lambda $ satisfying $\lambda
^{-1}n^{\frac{1+3\tau}{2}}p ^{3/2}=o(n^{1-\kappa})$ can be used.
Note that $\kappa+\tau/2<1/2$ by assumption. So in particular, we
can choose any $\lambda $ satisfying $\lambda (p
^{3/2}n)^{-1}\rightarrow\infty$ as $n\rightarrow\infty$.

\subsection{Proof of Theorem 3}
Theorem 3 is a straightforward corollary to Theorem 2 by the
argument in Step 2 of the proof of Theorem 1.

\bigskip

Throughout Sections A.4--A.6 below, we assume that $p >n$ and the
distribution of $\bz$ is continuous and spherically symmetric, that
is, invariant under the orthogonal group $\mathcal{O}(p )$. For
brevity, we use $\mathscr{L}\left(\cdot\right)$ to denote the
probability law or distribution of the random variable indicated.
Let $S^{q-1}(r)=\{x\in\mathbf{R}^q:\|x\|=r\}$ be the centered sphere
with radius $r$ in $q$-dimensional Euclidean space $\mathbf{R}^q$.
In particular, $S^{q-1}$ is referred to as the unit sphere in
$\mathbf{R}^q$.

\subsection{The distribution of $\bS = \left(\bZ\t\bZ\right)^+\bZ\t\bZ$}
It is a classical fact that the orthogonal group $\mathcal{O}(p )$
is compact and admits a probability measure that is invariant under
the action of itself, say,
\[ Q\cdot g\heq Q g,\quad g\in\mathcal{O}(p ),Q\in\mathcal{O}(p ).  \]
This invariant distribution is referred to as the uniform
distribution on the orthogonal group $\mathcal{O}(p )$. We often
encounter projection matrices in multivariate statistical analysis.
In fact, the set of all $p \times p $ projection matrices of rank
$n$ can equivalently be regarded as the Grassmann manifold
$\mathcal{G}_{p ,n}$ of all $n$-dimensional subspaces of the
Euclidean space $\mathbf{R}^{p }$; throughout, we do not distinguish
them and write
\[
\mathcal{G}_{p ,n}=\left\{U\t\diag\left(I_n
,0\right)U:U\in\mathcal{O}(p )\right\}.
\]

It is well known that the Grassmann manifold $\mathcal{G}_{p ,n}$ is
compact and there is a natural $\mathcal{O}(p )$-action on it, say,
\[ Q\cdot g\heq Q\t g Q,\quad g\in\mathcal{G}_{p ,n},Q\in\mathcal{O}(p ).  \]
Clearly, this group action is transitive, i.e. for any
$g_1,g_2\in\mathcal{G}_{p ,n}$, there exists some $Q\in\mathcal{O}(p
)$ such that $Q\cdot g_1=g_2$.  Moreover, $\mathcal{G}_{p ,n}$
admits a probability measure that is invariant under the
$\mathcal{O}(p )$-action defined above. This invariant distribution
is referred to as the uniform distribution on the Grassmann manifold
$\mathcal{G}_{p ,n}$.  For more on group action and invariant
measures on special manifolds, see Eaton (1989) and Chikuse (2003).

The uniform distribution on the Grassmann manifold is not easy to
deal with directly.  A useful fact is that the uniform distribution
on $\mathcal{G}_{p ,n}$ is the image measure of the uniform
distribution on $\mathcal{O}(p )$ under the mapping
\[ \varphi:\mathcal{O}(p )\rightarrow \mathcal{G}_{p ,n},\quad \varphi(U)=
U\t\diag\left(I_n ,0\right)U,\ U\in\mathcal{O}(p ).  \]

By the assumption that $\bz$ has a continuous distribution, we can
easily see that with probability one, the $n\times p $ matrix $\bZ$
has full rank $n$.  Let $\sqrt{\mu_1},\cdots,\sqrt{\mu_n }$ be its
$n$ singular values.  Then, $\bZ$ admits a singular value
decomposition
\begin{equation} \label{087}
\bZ=\bV\bD_1\bU,
\end{equation}
where $\bV\in\mathcal{O}(n)$, $\bU\in\mathcal{O}(p )$, and $\bD_1$
is an $n\times p $ diagonal matrix whose diagonal elements are
$\sqrt{\mu_1},\cdots,\sqrt{\mu_n }$, respectively.  Thus,
\begin{equation} \label{088}
\bZ\t \bZ=\bU\t\diag\left(\mu_1,\cdots,\mu_n ,0,\cdots,0\right)\bU
\end{equation}
and its Moore-Penrose generalized inverse is \[ \left(\bZ\t
\bZ\right)^+ =\sum_{i=1}^n\frac{1}{\mu_i}\bu_i\bu_i\t,
\]
where $\bU\t=\left(\bu_1,\cdots,\bu_{p }\right)$.  Therefore, we have
the following decomposition,
\begin{equation} \label{015}
\bS=\left(\bZ\t \bZ\right)^+\bZ\t \bZ=\bU\t\diag\left(I_n
,0\right)\bU,\quad \bU\in\mathcal{O}(p ).
\end{equation}

From (\ref{087}), we know that
$\bZ=\bV\diag\left(\sqrt{\mu_1},\cdots,\sqrt{\mu_n }\right)\left(I_n
,\bzero\right)_{n\times p }\bU$, and thus \[ \left(I_n
,\bzero\right)_{n\times p
}\bU=\diag\left(1/\sqrt{\mu_1},\cdots,1/\sqrt{\mu_n
}\right)\bV\t\bZ.
\] By the assumption that
$\mathscr{L}\left(\bz\right)$ is invariant under the orthogonal
group $\mathcal{O}(p )$, the distribution of $\bZ$ is also invariant
under $\mathcal{O}(p )$, i.e.,
\[
\bZ Q\deq \bZ\quad\text{for any }Q\in\mathcal{O}(p ).
\]
Thus, conditional on $\bV$ and $\left(\mu_1,\cdots,\mu_n \right)\t$,
the conditional distribution of $\left(I_n ,\bzero\right)_{n\times p
}\bU$ is invariant under $\mathcal{O}(p )$, which entails that
\[ \left(I_n ,\bzero\right)_{n\times p }\bU\deq\left(I_n ,\bzero
\right)_{n\times p }\widetilde{\bU}, \]
where $\widetilde{\bU}$ is uniformly distributed on the orthogonal
group $\mathcal{O}(p )$. In particular, we see that
$\left(\mu_1,\cdots,\mu \right)\t$ is independent of $\left(I_n
,\bzero\right)_{n\times p }\bU$. Therefore, these facts along with
(\ref{015}) yield the following lemma.

\begin{lemma}
$\mathscr{L}\left(\left(I_n ,\bzero\right)_{n\times p }\bU\right)
=\mathscr{L}\left(\left(I_n ,\bzero\right)_{n\times p
}\widetilde{\bU}\right)$ and $\left(\mu_1,\cdots,\mu_n \right)\t$ is
independent of $\left(I_n ,\bzero\right)_{n\times p }\bU$, where
$\widetilde{\bU}$ is uniformly distributed on the orthogonal group
$\mathcal{O}(p )$ and $\mu_1,\cdots,\mu_n $ are $n$ eigenvalues of
$\bZ\bZ\t$. Moreover, $\bS$ is uniformly distributed on the
Grassmann manifold $\mathcal{G}_{p ,n}$.
\end{lemma}

For simplicity, we do not distinguish $\widetilde{\bU}$ and $\bU$ in
the above singular value decomposition (\ref{087}).

\subsection{Deviation inequality on $\left\langle\bS\be_1,\be_1\right\rangle$}
\begin{lemma}
$\mathscr{L}\left(\left\langle\bS\be_1,\be_1\right\rangle\right)
=\frac{\chi_n^2}{\chi_n^2+\chi_{p -n}^2}$, where $\chi_n^2$ and
$\chi_{p -n}^2$ are two independent $\chi^2$-distributed random
variables with degrees of freedom $n$ and $p -n$, respectively, that
is, $\left\langle\bS\be_1,\be_1\right\rangle$ has a beta
distribution with parameters $n/2$ and $\left(p -n\right)/2$.
\end{lemma}

\begin{proof}
Lemma 1 gives $\mathscr{L}\left(\bS\right)=
\mathscr{L}\left(\bU\t\diag\left(I_n,0\right)\bU\right)$, where
$\bU$ is uniformly distributed on $\mathcal{O}(p )$.  Clearly,
$\left(\bU\be_1\right)$ is a random vector on the unit sphere $S^{p
-1}$.   It can be shown that $\bU\be_1$ is uniformly distributed on
the unit sphere $S^{p -1}$.

Let $\bW = (W_1,\cdots,W_{p })^T \sim {\cal N}(\bzero, I_p)$. Then,
we have $ \bU\be_1\deq  \bW / \| \bW\|$ and
\[
\left\langle\bS\be_1,\be_1\right\rangle
=\left(\bU\be_1\right)\t\diag\left(I_n,0\right)\bU\be_1 \deq
\frac{W_1^2+\cdots+V_n^2}{W_1^2+\cdots+W_{p }^2}.
\]
This proves Lemma 2.
\end{proof}

Lemmas 3 and 4 below give sharp deviation bounds on the
beta-distribution.

\begin{lemma}
{\em (Moderate deviation).} Let $\xi_1,\cdots,\xi_n$ be i.i.d.
$\chi_1^2$-distributed random variables. Then,
\begin{itemize}
\item[(i)] for any $\varepsilon>0$, we have
\begin{equation} \label{004}
P\left(n^{-1} ({\xi_1+\cdots+\xi_n})>1+\varepsilon\right)\leq
e^{-A_\varepsilon n},
\end{equation}
where $A_\varepsilon=\left[\varepsilon - \log ( 1 + \varepsilon)
\right]/2>0$.

\item[(ii)] for any $\varepsilon\in\left(0,1\right)$, we have
\begin{equation} \label{005}
P\left(n^{-1} ({\xi_1+\cdots+\xi_n}) <1-\varepsilon\right)\leq
e^{-B_\varepsilon n},
\end{equation}
where $B_\varepsilon=\left[-\varepsilon  - \log (1 - \varepsilon)
\right]/2>0$.
\end{itemize}
\end{lemma}

\begin{proof}
(i) Recall that the moment generating function of a
$\chi_1^2$-distributed random variable $\xi$ is
\begin{equation} \label{006}
M(t)=Ee^{t\xi}=\left(1-2t\right)^{-1/2},\quad
t\in\left(-\infty,1/2\right).
\end{equation}
Thus, for any $\varepsilon>0$ and $0<t<1/2$, by Chebyshev's
inequality (see, e.g. van der Vaart and Wellner, 1996) we have
\begin{align*}
P\left(\frac{\xi_1+\cdots+\xi_n}{n}>1+\varepsilon\right)
&\leq\frac{1}{e^{(t+1) n\varepsilon}}E\exp\left\{t\left(\xi_1
+\cdots+ \xi_n\right)\right\} = \exp(- n f_\varepsilon(t) ),
\end{align*}
where
$f_\varepsilon\left(t\right)=\frac{1}{2}\log\left(1-2t\right)+\left(1+\varepsilon\right)t$.
Setting the derivative $f_\varepsilon'\left(t\right)$ to zero gives
$t=\frac{\varepsilon}{2\left(1+\varepsilon\right)}$, where
$f_\varepsilon$ attains the maximum $A_\varepsilon=\left[\varepsilon
- \log(1+\varepsilon) \right]/2$, $\varepsilon>0$. Therefore, we
have
\[ P\left(n^{-1} ({\xi_1+\cdots+\xi_n})>1+\varepsilon\right)\leq
e^{-A_\varepsilon n}.  \] This proves (\ref{004}).

(ii) For any $0<\varepsilon<1$ and $t>0$, by Chebyshev's inequality
and (\ref{006}), we have
\begin{align*}
P\left(n^{-1} ({\xi_1+\cdots+\xi_n}) <1-\varepsilon\right)
&\leq\frac{1}{e^{tn\varepsilon}}E\exp\left\{t\left(1-\xi_1\right)+\cdots+
t\left(1-\xi_n\right)\right\} = \exp(-ng_\varepsilon(t)),
\end{align*}
where
$g_\varepsilon\left(t\right)=\frac{1}{2}\log\left(1+2t\right)-\left(1-\varepsilon\right)t$.
Taking $t=\varepsilon/ ({2\left(1-\varepsilon\right)})$ yields
(\ref{005}).
\end{proof}

\begin{lemma}
{\em (Moderate deviation).} For any $C>0$, there exist constants
$c_1$ and $c_2$ with $0<c_1<1<c_2$ such that
\begin{equation} \label{008}
P\left(\left\langle\bS\be_1,\be_1\right\rangle<c_1\frac{n}{p }
\text{ or }>c_2\frac{n}{p }\right)\leq 4e^{-Cn}.
\end{equation}
\end{lemma}

\begin{proof}
From Lemma 3, we know that $
\left\langle\bS\be_1,\be_1\right\rangle\deq {\xi}/ {\eta}, $ where
$\xi$ is $\chi_n^2$-distributed and $\eta$ is $\chi_{p
}^2$-distributed.  Note that $A_\varepsilon$ and $B_\varepsilon$ are
increasing in $\varepsilon$ and have the same range $(0,\infty)$.
For any $C>0$, it follows from the proof of Lemma 3 that there exist
$\widetilde{c}_1$ and $\widetilde{c}_2$ with
$0<\widetilde{c}_1<1<\widetilde{c}_2$, such that
$B_{1-\widetilde{c}_1}=C$ and $A_{\widetilde{c}_2-1}=C$.  Now define
\[ \mathcal{A}=\left\{\frac{\xi}{n}<\widetilde{c}_1\text{ or
}>\widetilde{c}_2\right\}\quad \text{and}\quad
\mathcal{B}=\left\{\frac{\eta}{p }<\widetilde{c}_1\text{ or
}>\widetilde{c}_2\right\}.  \] Let
$c_1=\widetilde{c}_1/\widetilde{c}_2$ and
$c_2=\widetilde{c}_2/\widetilde{c}_1$. Then, it can easily be shown
that
\begin{equation} \label{010}
\left\{\left\langle\bS\be_1,\be_1\right\rangle<c_1\frac{n}{p }
\text{ or }>c_2\frac{n}{p
}\right\}\subset\mathcal{A}\cup\mathcal{B}.
\end{equation}
It follows from (\ref{004}) and (\ref{005}) and the choice of
$\widetilde{c}_1$ and $\widetilde{c}_2$ above that
\begin{equation} \label{016}
P\left(\mathcal{A}\right)\leq2e^{-Cn}\quad\text{and}\quad
P\left(\mathcal{B}\right)\leq2e^{-Cp }.
\end{equation}
Therefore, by $p \geq n$ and Bonferroni's inequality, the results
follow from (\ref{010}) and (\ref{016}).
\end{proof}

\subsection{Deviation inequality on $\left\langle\bS\be_1,\be_2\right\rangle$}
\begin{lemma}
Let $\bS\be_1=\left(V_1,V_2,\cdots,V_{p }\right)\t$. Then, given
that the first coordinate $V_1=v$, the random vector
$\left(V_2,\cdots,V_{p }\right)\t$ is uniformly distributed on the
sphere $S^{p -2}(\sqrt{v-v^2})$.  Moreover, for any $C>0$, there
exists some $c>1$ such that
\begin{equation} \label{020}
P\left(\left|V_2\right|> c \sqrt{n}p
^{-1}\left|W\right|\right)\leq3e^{-Cn},
\end{equation}
where $W$ is an independent $\mathcal{N}(0,1)$-distributed random
variable.
\end{lemma}

\begin{proof} In view of (\ref{015}),  it follows that
\begin{align*}
\|bV\|^2  =\be_1\t\bS\be_1 =V_1,
\end{align*}
where $\bV = \left(V_1,\cdots,V_{p }\right)\t$.  For any
$Q\in\mathcal{O}(p -1)$, let $\widetilde{Q}=\diag\left(1,Q\right)
\in\mathcal{O}(p )$. Thus, by Lemma 1, we have
\begin{align*}
\widetilde{Q} \bV  &\deq
=\left(\bU\widetilde{Q}\t\right)\t\diag\left(I_n,0\right)\left(
\bU\widetilde{Q}\t\right)\widetilde{Q}\be_1\\
&\deq\bU\t\diag\left(I_n,0\right)\bU\be_1\deq \bV.
\end{align*}
This shows that given $V_1=v$, the conditional distribution of
$\left(V_2,\cdots,V_{p }\right)\t$ is invariant under the orthogonal
group $\mathcal{O}(p -1)$.   Therefore, given $V_1=v$, the random
vector $\left(V_2,\cdots,V_{p }\right)\t$ is uniformly distributed
on the sphere $S^{p -2}(\sqrt{v-v^2})$.

Let $W_1,\cdots,W_{p -1}$ be i.i.d.  $\mathcal{N}(0,1)$-distributed
random variables, independent of $V_1$.  Conditioning on $V_1$, we
have
\begin{equation} \label{022}
V_2\deq\sqrt{V_1-V_1^2}\frac{W_1}{\sqrt{W_1^2+\cdots+W_{p -1}^2}}.
\end{equation}
Let $C>0$ be a constant.  From the proof of Lemma 4, we know that
there exists some $c_2>1$ such that
\begin{equation} \label{018}
P\left(V_1>c_2{n}/{p }\right)\leq2e^{-Cn}.
\end{equation}
It follows from (\ref{005}) that there exists some $0<c_1<1$ such
that
\begin{equation} \label{019}
P\left({W_1^2+\cdots+W_{p -1}^2}<c_1({p -1})\right)\leq e^{-C\left(p
-1\right)}\leq e^{-Cn},
\end{equation}
since $p >n$.  Let $c=\sqrt{c_2/c_1}$. Then, by $V_1-V_1^2\leq V_1$
and Bonferroni's inequality, (\ref{020}) follows immediately from
(\ref{022})--(\ref{019}).
\end{proof}

\subsection{Verifying Property C for Gaussian distributions}
In this section, we check Property C in (\ref{086}) for Gaussian
distributions. Assume $\bx$ has a $p $-variate Gaussian
distribution. Then, the $n\times p $ design matrix
$\bX\sim\mathcal{N}(\bzero,I_n\otimes\Sig )$ and
\[ \bZ\sim\mathcal{N}(\bzero,I_n\otimes
I_{p })=\mathcal{N}(\bzero,I_{n\times p }), \] i.e., all the entries
of $\bZ$ are i.i.d. $\mathcal{N}(0,1)$ random variables, where the
symbol $\otimes$ denotes the Kronecker product of two matrices. We
will invoke results in the random matrix theory on extreme
eigenvalues of random matrices in Gaussian ensemble.

Before proceeding, let us make two simple observations. First, in
studying singular values of $\bZ$, the role of $n$ and $p $ is
symmetric. Second, when $p >n$, by letting
$\bW\sim\mathcal{N}(\bzero,I_{m\times p })$, independent of $\bZ$,
and
\[ \widetilde{\bZ}_{\left(n+m\right)\times p }=\left(\begin{array}{c}
\bZ \\
\bW
\end{array}\right), \]
then the extreme singular values of $\bZ$ are sandwiched by those of
$\widetilde{\bZ}$. Therefore, a combination of Lemmas 6 and 7 below
immediately implies Property C in (\ref{086}).

\begin{lemma}
Let $p \geq n$ and $\bZ\sim\mathcal{N}(\bzero,I_{n\times p })$.
Then, there exists some $C>0$ such that for any eigenvalue $\lambda$
of $p ^{-1}\bZ\bZ\t$ and any $r>0$,
\[ P\left(\left|\sqrt{\lambda}-E(\sqrt{\lambda})\right|>r\right)\leq Ce^{-p r^2/C}. \]
Moreover, for each $\lambda$, the same inequality holds for a median
of $\sqrt{\lambda}$ instead of the mean.
\end{lemma}

\begin{proof}
See Proposition 3.2 in Ledoux (2005) and note that Gaussian measures
satisfy the dimension-free concentration inequality (3.6) in Ledoux
(2005).
\end{proof}

\begin{lemma}
Let $\bZ\sim\mathcal{N}(\bzero,I_{n\times p })$. If $p
/n\rightarrow\gamma>1$ as $n\rightarrow\infty$, then we have
\[
\lim_{n\rightarrow\infty}\text{median}\left(\sqrt{\lambda_{\text{max}}(p
^{-1}\bZ\bZ\t)}\right)= 1+\gamma^{-1/2} \] and
\[
\liminf_{n\rightarrow\infty}E\left(\sqrt{\lambda_{\text{min}}(p
^{-1}\bZ\bZ\t)}\right)\geq 1-\gamma^{-1/2}.
\]
\end{lemma}

\begin{proof}
The first result follows directly from Geman (1980):
\[ \lambda_{\text{max}}(p ^{-1}\bZ\bZ\t)\toas
\left(1+\gamma^{-1/2}\right)^2\quad\text{as }n\rightarrow\infty. \]
For the smallest eigenvalue, it is well known that (see, e.g.,
Silverstein, 1985 or Bai, 1999)
\[ \lambda_{\text{min}}(p ^{-1}\bZ\bZ\t)\toas
\left(1-\gamma^{-1/2}\right)^2\quad\text{as }n\rightarrow\infty. \]
This and  Fatou's lemma entails the second result.
\end{proof}

\end{document}